\documentclass[12pt]{article}

\usepackage{amsmath,amssymb}

\usepackage[
]{hyperref}
\hypersetup{
	colorlinks,
	linkcolor={red!80!black},
	citecolor={blue!50!black},
	urlcolor={blue!80!black}
}
\usepackage{xcolor}

\usepackage{easyReview}

\usepackage{fullpage}

\usepackage{enumitem}

\usepackage{bbold}

\usepackage{xcolor,graphicx}
\usepackage{wrapfig}
\usepackage{caption}
\usepackage{subcaption}
\usepackage{cleveref}
\usepackage{amsthm}
\theoremstyle{definition}
\newtheorem{theorem}{Theorem}[section]
\newtheorem{proposition}[theorem]{Proposition}
\newtheorem{lemma}[theorem]{Lemma}

\newtheorem{conjecture}[theorem]{Conjecture}
\theoremstyle{remark}
\newtheorem{remark}{Remark}[section]
\newtheorem{example}{Example}

\usepackage[tworuled, lined, linesnumbered]{algorithm2e}
\SetKwFor{While}{while}{\emph{:}}{}%
\SetKwFor{For}{for}{\emph{\string:}}{}%
\SetKwFor{ForEach}{for each}{\emph{\string:}}{}%
\SetKwIF{If}{ElseIf}{Else}{if}{\emph{:}}{else if}{else\emph{:}}{}%
\SetStartEndCondition{ }{}{}
\DontPrintSemicolon
\SetVlineSkip{0pt}
\SetAlgoInsideSkip{medskip}
\SetAlgoSkip{bigskip}
\SetKwSty{emph}

\usepackage{etoolbox}
\makeatletter
\patchcmd{\@algocf@start}
{-1.5em}
{0pt}
{}{}
\makeatother

\usepackage{authblk}

\usepackage{fontawesome5}

\begin{document}
\title{A dynamical approach to spanning and surplus edges of random graphs}
\date{}
\author[1,2]{Josué Corujo}
\author[2]{Vlada Limic}

\affil[1]{
Univ Paris Est Créteil, Univ Gustave Eiffel, CNRS, LAMA UMR 8050, F-94010 Créteil, France
}
\affil[2]{Institut de Recherche Mathématique Avancée, UMR 7501 Université de Strasbourg et CNRS, 
7 rue René-Descartes, 67000 Strasbourg, France}
{
	\makeatletter
	\renewcommand\AB@affilsepx{: \protect\Affilfont}
	\makeatother
	
	\affil[ ]{Email ids}
	
	\makeatletter
	\renewcommand\AB@affilsepx{, \protect\Affilfont}
	\makeatother
	
	\affil[1]{\href{mailto:josue.corujo-rodriguez@u-pec.fr}{josue.corujo-rodriguez@u-pec.fr}}
	\affil[2]{\href{mailto:vlada@math.unistra.fr}{vlada@math.unistra.fr}}
}

\maketitle
\begin{abstract}
Consider a finite inhomogeneous random graph evolving in continuous time, where each vertex is assigned a mass, and an edge between any pair of vertices appears at a rate proportional to the product of their masses. The process tracking the evolution of component sizes evolves according to the multiplicative coalescent dynamic and can be encoded using the simultaneous breadth-first walk introduced by Limic (2019).
We extend this encoding to incorporate surplus edge data within each connected component. Two distinct graph-based representations of the multiplicative coalescent, each with its own advantages and limitations, are analyzed in detail. In particular, a canonical multigraph introduced by Bhamidi, Budhiraja and Wang (2014), which is naturally connected to the augmented multiplicative coalescent, emerges from our framework.
We demonstrate that a transformation of the simultaneous breadth-first walk, supplemented with an additional and independent source of randomness, encodes the full dynamics of the augmented multiplicative coalescent. 
\end{abstract}





\section{Introduction}

For $n \in \mathbb{N}$, write $[n]$ for $\{1,\ldots, n\}$.
Let us denote by $G(n,p)$, where $p\in[0,1]$, the  (a.k.a.\ Erd\H{o}s\,--\,Rényi\,--\,Gilbert)
random graph, originally introduced by Gilbert \cite{1959Gilbert} and first studied by \cite{erdren}, where each edge is included in the graph with probability $p$, independently from every other edge.
A continuous-time variation of $G(n,p)$, where $p\in[0,1]$, is naturally constructed as follows: 
fix the $n$ vertices $[n]$ and let each of the $\binom{n}{2}$ edges appear at an exponential time of rate $1$, independently of each other.
This graph-valued Markov chain was first introduced by Stepanov \cite{Stepanov1970OnTP} and we call it Erd\H{o}s\,--\,Rényi\,--\,Stepanov random graph process.
This transforms the model into a continuous-time Markov chain, running on the set of graphs with vertices $[n]$ and going from the trivial graph ($n$ disconnected vertices) at time $t = 0$, to the complete graph when $t \to \infty$.
This continuous-time construction is equally obtained by the time-change $t = -\ln(1-p)$ in the natural coupling of $\big(G(n,p), p \in [0,1]\big)$, i.e.\ the process where at each $p \in [0,1]$ the edge $\{i,j\}$ is present in the graph if and only if $\{ U_{\{i,j\}} \le p \}$, for a fixed family $ \displaystyle ( U_{\{i,j\}} )_{i < j}$ of independent random variables with uniform distribution in $[0,1]$.
 
Here and in the rest of the paper {\em connected} means connected by a path of edges in the usual graph theory sense. 
If the minimal path is in fact an edge, this is typically clear from the context. A connected component is a subset $S$ of vertices such that any two vertices in $S$ are connected, and no vertex in $S^c$ is connected to a vertex in $S$.
With this convention, any two different connected components merge at the minimal connection time of a pair of vertices $(k,l)$ (where $k$ is from one, and $l$ from the other component) to form a single connected component.
Let the mass of any connected component be equal to the number of its vertices.
Due to elementary properties of independent exponential random variables, it is immediate that 
\begin{equation}
\label{merge}
\begin{array}{c}
\mbox{ any pair of components with masses (sizes) $x$ and $y$ merges }\\
\mbox{ at rate $x \cdot y$ into a single component of mass $x+y$.}
\end{array}
\end{equation}
In other words, the vector of ordered sizes of the connected components of the  continuous-time random graph evolves according to the \emph{multiplicative coalescent} (MC) dynamics.
Due to the relation $p= 1-\mathrm{e}^{-t}$ (and $1 - \mathrm{e}^{-t} \approx t$ for small $t$), this continuous-time random graph exhibits the same phase transition as $G(n,\cdot)$ does as $n$ diverges, at the critical time $1/n$ plus a lower order term. 

Aldous \cite{aldRGMC} extended this construction as follows: instead of mass $1$, let vertex $i\in [n]$ have initial mass $x_i>0$. 
For each $i,j\in [n]$ let the edge between $i$ and $j$ appear at rate $x_i \cdot x_j$, independently of others.
In the sequel, we will call this process the \emph{Aldous' (inhomogeneous) continuous-time random graph} and denote by $(\mathcal{G}^{\pmb{x}}(t), t\ge 0)$ a graph-valued continuous-time Markov process following this dynamic.

A closely related random graph model (in fact, it is a time-changed version of our model), was called \emph{multiplicative graph} by Broutin et al.\ \cite{Broutin2021,Broutin2022}, and analyzed by the same authors using a depth-first exploration process. 
The breadth-first search we use in this paper allows a dynamical study of the random graph process. 
As far as we know, there is no evidence that any depth-first algorithm is capable of doing such analysis. 
However, the depth-first walk is a stronger tool for understanding the topological properties (in particular, the scaling limit of the endowed distance) of a near critical random graph at a fixed (near-critical) time.

Notice that the transition mechanism of merging in the Aldous' continuous-time random graph is again \eqref{merge}.
Furthermore, Proposition 4 in \cite{aldRGMC}  shows that if the set of vertices is ${\mathbb N}$, and if $\pmb{x}=(x_1,x_2,\ldots) \in l^2$, where the initial mass of $i$ is $x_i$, then this (infinite) random graph process is still well-defined, and its connected component masses form a vector in $l^2$ at any later time a.s. 
Here is a more precise formulation: let $(l^2_\searrow,d)$ be the metric space of infinite sequences
${\pmb{y}} = (y_1,y_2,\ldots)$
such that $y_1 \geq y_2 \geq \ldots \geq 0$
and $\sum_i y_i^2 < \infty$,
with the distance 
$\operatorname{d}({\pmb{y}},\pmb{ z} ) = \sqrt{\sum_i (y_i-z_i)^2}$.
Let ``$\operatorname{ord}$'' be the ``decreasing
ordering'' map defined on infinite-length vectors, and note that it is well defined on $l^2_\searrow$.
Let $X_i(t)$ be the $i^{\text{th}}$ largest connected component mass in $\mathcal{G}^{\pmb{x}}(t)$, for every $t \ge 0$.   
The process $(\pmb{X}(t),\,t\geq 0)\equiv \big((X_1(t),X_2(t),\ldots),t\geq 0\big)$ started from $\pmb{X}(0) = {\rm ord}(\pmb{x}) \in l^2_\searrow$ is an $l^2_\searrow$-valued Feller process evolving according to \eqref{merge}.
See \cite[Prop.\ 4 and  5]{aldRGMC}, and \cite[\S\ 2.1]{VLthesis} for an alternative derivation of the Feller property.
Starting with \cite{aldRGMC}, any such process $\pmb{X}$ is referred to as a {\em multiplicative coalescent}.
In this note, a {\em graph representation of the multiplicative coalescent} (or an {\em MC graph representation} for short) will be any random graph-valued process such that its corresponding ordered component sizes evolve as a multiplicative coalescent.

Aldous' continuous-time random graph $(\mathcal{G}^{\pmb{x}}(t), t\ge 0)$ is clearly a MC graph representation.
A different but similar MC graph, which we will now recall,  was explored by Bhamidi et al.~in \cite[\S\ 2.3.1]{bhamidietal2}.
Here for each $i,j\in [n]$ a new {\em directed edge}
 $i \rightarrow j$ appears at rate $x_i \cdot x_j/2$, and for each $i$ a self-loop $i \rightarrow i$ appears at rate $x_i^2/2$. 
This random-graph is in fact an oriented multi-graph, since it is an oriented graph with loops and multiple edges allowed.
Let us denote by $(\mathcal{MG}^{\pmb{x}}(t), t\ge 0)$, a continuous-time multi-graph-valued Markov process following this dynamic.
If the connected components are obtained by taking into account all the edges (regardless of their orientation), and the mass of each connected component is again the sum of masses of its participating vertices,
it is easy to see that the resulting ordered component masses evolve again according to the {multiplicative coalescent} transitions. 
Indeed, the presence of multi-edges and loops does not change the connectivity properties or the component masses, so the random graph process from Bhamidi et al.~\cite{bhamidietal2} can be  coupled to Aldous' construction outlined in the previous paragraph.
Furthermore, one could also record the information on the surplus edges of the connected components in $(\mathcal{MG}^{\pmb{x}}(t), t\ge 0)$.
Let
\(
	\mathbb{N}^\infty = \{ (n_1, n_2, \dots): n_i \in \mathbb{N}_0, \text{ for all } i \ge 1 \},
\)
and define
\begin{equation*} 
	\mathbb{U}_\searrow = \left\{ \big( (x_i)_{i \ge 1}, (n_i)_{i \ge 1} \big) \in l^2_\searrow \times \mathbb{N}^{\infty} : \text{ if } x_i = 0, \, \text{ then } n_i = 0, \text{ and  if } x_i = x_{i+1}, \text{ then } n_i \ge n_{i + 1} \right\}.
\end{equation*}
For every element $(x_i, n_i)_{i \ge 1} \in \mathbb{U}_\searrow$, the positive real $x_i$ represents the mass of the $i^{\mathrm{th}}$ largest component, while the positive integer $n_i$ is the number of surplus edges in the same component.
In this paper we consider only the case where the elements in $\mathbb{U}_\searrow$ are finite, i.e.\ they only have a finite number of non zero components.
The definition of $\mathbb{U}_\searrow$ when considering infinite length vectors is more subtle (see \cite{bhamidietal2} for more details), however not important for the present work.

Hence, the joint evolution of component
masses and surplus edge counts in $(\mathcal{MG}^{\pmb{x}}(t), t\ge 0)$ has the following transitions:
\begin{description}
 \item [coalescence jump:] for each $i \neq j$,
 at rate $x_i \cdot x_j$,
the process jumps from $(\pmb{x}, \pmb{n})$ to 
 $(\pmb{x}^{i, j}, \pmb{n}^{i,j})$ 
 where $(\pmb{x}^{i, j}, \pmb{n}^{i,j})$  is obtained by merging components $i$ and $j$ into a component with mass $x_i + x_j$ and surplus $n_i + n_j$, followed by reordering the coordinates with respect to the masses in order to obtain again an element of $\mathbb{U}_\searrow$,
 \item [surplus jump:] for each $i \ge 1$, at rate $x_i^2/2$,
 the process jumps from $(\pmb{x}, \pmb{n})$ to 
 $(\pmb{x}, \pmb{n}^{i})$, where $(\pmb{x}, \pmb{n}^{i})$ is the state obtained by changing only the $i^{\mathrm{th}}$ component $(x_i,n_i)$ of $(\pmb{x}, \pmb{n})$ into $(x_i,n_i+1)$,
and reordering the coordinates if needed, to obtain an element in $\mathbb{U}_\searrow$.
\end{description}
We will call a process with this dynamic an \emph{augmented multiplicative coalescent}.
See \cite{bhamidietal2} for a detailed study of these processes.

As already hinted above,  in our setting it is convenient to embed finite vectors into an infinite-dimensional space.
Refer henceforth to $\pmb{x}=(x_1,x_2,x_3\ldots)\in l^2_\searrow$ as {\em finite}, if  for some $i\in \mathbb{N}$ we have $x_i=0$. Let the {\em length} of $\pmb{x}$ be the number $\operatorname{len}(\pmb{x})$ of non-zero coordinates of $\pmb{x}$.
We denote by $\pmb{0}$ the infinite vector of zeros.
Fix a finite initial configuration $\pmb{x}\in l^2_\searrow$.
For each $i\leq \operatorname{len}(\pmb{x})$, let $\xi_i$ have exponential distribution with rate $x_i$, independently over $i$. 

The order statistics of $(\xi_i)_{i\leq \operatorname{len}(\pmb{x})}$ are denoted by $(\xi_{(i)})_{i\leq \operatorname{len}(\pmb{x})}$, and let us denote by $(\pi_i)_{i \leq \operatorname{len}(\pmb{x})}$ the permutation induced by the ordering of $(\xi_{(i)})_{i\leq \operatorname{len}(\pmb{x})}$, i.e.\ 
\[
	\pi_i = k \mbox{ if and only if } \xi_{(i)} = \xi_{k}.
\]
In this way, $(x_{\pi_1},x_{\pi_2},\ldots, x_{\pi_{\operatorname{len}(\pmb{x})} } )$  is the size-biased random ordering of the initial non-trivial block masses.
Given $\pmb{\xi} := (\xi_1, \dots, \xi_{\operatorname{len}(\pmb{x})})$, define simultaneously, for all $s \ge 0$ and $q>0$, the process
\begin{align}
	Z^{\pmb{x},q}(s) :=& \sum_{i=1}^{\operatorname{len}(\pmb{x})} x_i \cdot \mathbb{1}_{(\xi_i/q\, \leq \ s)} -s \nonumber \\
	=& \sum_{i=1}^{\operatorname{len}(\pmb{x})} x_{\pi_i} \cdot \mathbb{1}_{(\xi_{(i)}/q\, \leq \ s)} -s. \label{defZbxq}
\end{align}
In words, $Z^{\pmb{x},q}$ has a unit negative drift and successive positive jumps, which occur precisely at times $(\xi_{(i)}/q)_{i\leq \operatorname{len}(\pmb{x})}$, and where the $i^{\text{th}}$  successive jump is of magnitude $x_{\pi_i}$.
Note that, as $q$ approaches zero, the exponential jump times $\xi_\cdot/q$ diverge, but more importantly the distances between any two successive jump times also diverge, which is consistent with the absence of edges in the random graph at small times.

In \cite{multcoalnew} the family $Z^{\pmb{x},\cdot}$ defined by \eqref{defZbxq} was called the {\em simultaneous breadth-first walks}. 
It was shown in \cite[Prop.\ 5]{multcoalnew} that, as $q$ increases, the excursion lengths of the reflected $Z^{\pmb{x},\cdot}$ have the law of the  {multiplicative coalescent} started from the configuration $\pmb{x}$. 
This was an essential step in the proof of  \cite[Thm.\ 1.2]{multcoalnew}, a characterization of the trajectories of the extreme eternal version of the MC. 
In that paper the family $Z^{\pmb{x},\cdot}$ was related to a graph representation of the {multiplicative coalescent}. 
Indeed, this MC graph representation recalled in Section \ref{S:Coupl} below
is a finite random forest whose tree masses evolve in time according to \eqref{merge}. 
Recall that for a finite connected graph $(V, E)$ with $|V| = n$, the number of surplus edges is defined as $|E| - (n-1)$, since any spanning tree of $(V,E)$ has exactly $n-1$ vertices.
The number of surplus edges are a measure of the level of connectivity of the connected components in the random graph.
Section \ref{S:Coupl} will introduce two extensions of the breadth-first search algorithm studied in \cite{multcoalnew} with the purpose of adding the surplus edges to each tree in the forest, relating the
MC forest-valued representation to the processes $(\mathcal{G}^{\pmb{x}}(q), q\ge 0)$ and $(\mathcal{MG}^{\pmb{x}}(q), q\ge 0)$.

We hereby wish to point out two crucial consequences arising from the research presented in this paper. 
In Section \ref{S:take2} we obtain 
a simultaneous breadth-first walk encoding of the marginal law of the augmented multiplicative coalescent at time $q > 0$, starting from a finite vector $\pmb{x}$ and no surplus edges.
Section \ref{sec:sBFWencodesAMCdynamic} is devoted to a modified encoding construction of the augmented multiplicative coalescent, as a process evolving in $q$.

Before stating these results, we introduce additional notation.
Define 
$$ B^{\pmb{x},q}(s):= Z^{\pmb{x},q}(s) - \inf_{u\leq s} Z^{\pmb{x},q}(u), \ s\geq 0,$$
and let $\Lambda$ be a homogeneous Poisson point process on $[0,\infty)\times [0,\infty)$, independent of $Z^{\pmb{x}, q}$.
Let $X_i^{\pmb{x}}(q)$ be the length of the $i^{\text{th}}$ largest excursion of $s \mapsto q \cdot B^{\pmb{x},q}(s)$ (or equivalently, of $B^{\pmb{x},q}$) above zero  and let $N_i^{\pmb{x}}(q)$ be the number of points in $\Lambda$ below the curve $s \mapsto q \cdot B^{\pmb{x},q}(s)$ during the above-mentioned excursion, such that 
\begin{equation}\label{eq:def_encoding}
    \big( \pmb{X}^{\pmb{x}}(q), \pmb{N}^{\pmb{x}}(q) \big) := \big( (X_1^{\pmb{x}}(q), X_2^{\pmb{x}}(q), \dots), (N_1^{\pmb{x}}(q), N_2^{\pmb{x}}(q), \dots) \big) \in \mathbb{U}_\searrow.
\end{equation}
Now we can state our first main result.
\begin{theorem}[Encoding of the marginal law of the AMC]
\label{thm:encodingAMC}
    The random vector $\big( \pmb{X}^{\pmb{x}}(q), \pmb{N}^{\pmb{x}}(q) \big)$ has the law of the augmented multiplicative coalescent evaluated at time $q > 0$, and started from $(\pmb{x},\pmb{0})$.
\end{theorem}
The proof of Theorem \ref{thm:encodingAMC} is postponed to Section \ref{S:Multi-g}.

Although Theorem  \ref{thm:encodingAMC} is a static result (in the sense that it characterizes the marginal laws of the AMC), its proof relies on the surplus construction given in Sections \ref{S:take2}, which is dynamical, meaning that it links the trajectories of the process $(\mathcal{MG}^{\pmb{x}}(t), t\ge 0)$ and the simultaneous breadth-first walk.
It is reasonable to expect that a modification to this encoding could encode the trajectories of the AMC.
This is indeed the case.

In order to explain this, we first introduce a 
modification of $B^{\pmb{x},q}$.
Notice that the curve $B^{\pmb{x},q}$ consists of excursions above zero separated by successive intervals where $B^{\pmb{x},q}$  remains at zero.
We call those intervals the {\em load-free intervals (or periods)}.
They are important part of the encoding in \cite{multcoalnew}: the specific gap between excursions 
which decreases in $q$ ensures
that the merging of excursions obeys the multiplicative coalescent dynamics.
Our new idea is to remove the load-free intervals $\omega$-by-$\omega$ in order to get a transformed curve, denoted $\mathfrak{B}^{\pmb{x}, q}$, and such that the time of merging of its excursions again obeys the multiplicative coalescent dynamics.
Note that outside the load-free intervals $Z^{\pmb{x},q}$ is larger than its past infimum, and moreover that
during the load-free intervals the past infimum $s \mapsto \inf_{u \le s} Z^{\pmb{x},q}(u)$ decreases at constant rate $1$.
Thus, if we denote by $\rho_s$ the cumulative length of the load-free periods of $B^{\pmb{x},q}$ up to time $s$, i.e.\
\begin{equation*} \label{eq:def_rho}
    \rho_s := \int_0^s \mathbb{1}_{\{B^{\pmb{x},q}(u) = 0\}} \mathrm{d}u,
\end{equation*}
we have the identity
\begin{equation*} \label{eq:charact_rho}
    \rho_s = - \inf_{u \le s} Z^{\pmb{x},q}(u).
\end{equation*}
We define $\mathfrak{B}^{\pmb{x}, q}$ as the process such that $\mathfrak{B}^{\pmb{x}, q}(s) = 0$, for every $s \ge \sum x_i$ and on $[0, \sum_i x_i]$ $\mathfrak{B}^{\pmb{x}, q}$ is defined via the following identities
\begin{equation} \label{eq:def_B_curly}
{B}^{\pmb{x},q}\left( s \right) = \mathfrak{B}^{\pmb{x},q}\left( s + \inf_{u \le s} Z^{\pmb{x},q}(s) \right) = \mathfrak{B}^{\pmb{x},q}\left( s -\rho_s \right) , \  \ s \ge 0, \, q > 0.
\end{equation}
In this way, the values of $\mathfrak{B}^{\pmb{x}, q}$ are completely determined, and the path of $\mathfrak{B}^{\pmb{x},q}$ is precisely the path of ${B}^{\pmb{x}, q}$ after all of its load-free intervals are removed. 
The process $\mathfrak{B}^{\pmb{x}, 0 }$, for $q = 0$, is defined as the limit of $\mathfrak{B}^{\pmb{x}, q}$ when $q\to 0$ and it simply consists on $\mathrm{len}(\pmb{x})$ triangle-shaped excursions such that the $i^{\text{th}}$ excursion has length $x_{\pi_i}$.

It is clear that $\mathfrak{B}^{\pmb{x},q}$ decreases at a constant unit drift on the entire interval $[0, \sum_{i \ge 1} x_i]$, and that the jump occurring at $\xi_i/q$ in the path of $B^{\pmb{x},q}$ is shifted to $\xi_i/q-\rho_{\xi_i/q}$ in the path of $\mathfrak{B}^{\pmb{x},q}$.
Hence, we have the following explicit and alternative expression for $\mathfrak{B}^{\pmb{x},q}$ on $\left[ 0, \sum_{i = 1}^n x_i \right]$
\begin{equation*}
\label{eq:def_B_curly2}
\mathfrak{B}^{\pmb{x},q}(s) = \sum_{i = 1}^n x_{i} \cdot \mathbb{1}_{\{ \Xi(\xi_i/q) \le s \}} - s, \;\; \text{ where }\ \Xi: s \mapsto s - \rho_{s}.   
\end{equation*}
Note that the function $\Xi$ can be also defined as
\begin{equation*}
\label{def:Xi}
    \Xi : s \mapsto \int_0^s \mathbb{1}_{\{ B^{\pmb{x}, q}(u) > 0 \}} \mathrm{d} u.
\end{equation*}
In particular, $\Xi$ is non-decreasing.

\begin{figure}[h]
    \centering
    \begin{subfigure}[b]{0.49\textwidth}
        \centering
        \includegraphics[width=\textwidth]{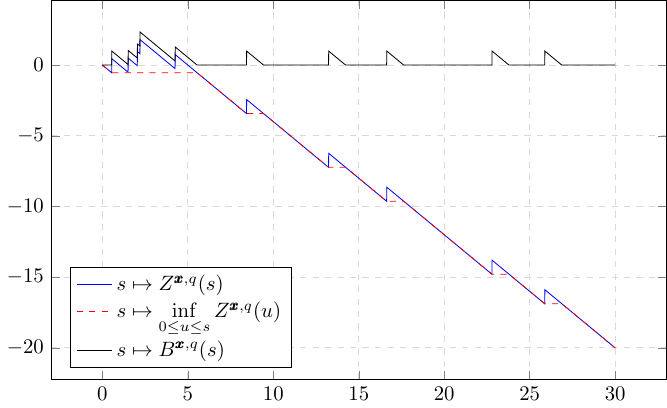}
        \caption{sBFW and its reflected version}
        \label{fig:transformed_B1}
    \end{subfigure}
    \hfill
    \begin{subfigure}[b]{0.49\textwidth}
        \centering
        \includegraphics[width=\textwidth]{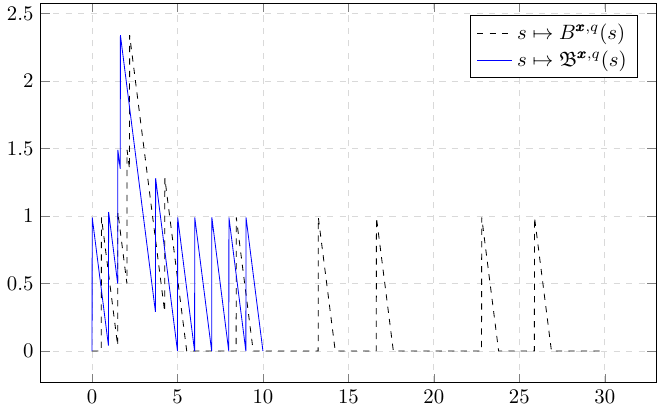}
        \caption{Reflected  sBFW and its modified version}
        \label{fig:transformed_B2}
    \end{subfigure}
    \caption{
    Plots of a simulation of the simultaneous breadth-first walk and its reflected version on the left panel, and the same reflected  simultaneous breadth-first walk together with its modified version obtained by removing the ``load-free intervals'', on the right panel. 
    The vector of initial masses $\pmb{x}$ is the $10$ dimensional all-one vector and $q = 0.7$.}
    \label{fig:modified_SBFW}
\end{figure}

Figures \ref{fig:transformed_B1} and \ref{fig:transformed_B2} show graphs of $Z^{\pmb{x},q}$, its reflected version $B^{\pmb{x},q}$ and its modification $\mathfrak{B}^{\pmb{x},q}$.
Notice how $\mathfrak{B}^{\pmb{x},q}$ is obtained from $B^{\pmb{x},q}$ by removing the load-free intervals. 

\begin{remark}[A link to the original BFW]
Note that the law of of $\mathfrak{B}^{\pmb{x}, q}$ is reminiscent of the law of the original BFW of \cite{aldRGMC,EBMC}, see also \cite[Lemma 2]{multcoalnew} (they differ in the presence/absence of jumps at the beginnings of excursions). 
\end{remark}

Let $\Lambda$ be a homogeneous Poisson point process on $[0,\infty)^2$ independent of $(\mathfrak{B}^{\pmb{x}, q},q\geq 0)$.
Define $\mathfrak{N}_i(q)$ as the number of points of $\Lambda$ in the region below the curve $s \mapsto q \cdot \mathfrak{B}^{\pmb{x}, q}(s)$ restricted to its $i^{\text{th}}$ longest excursion above $0$, whose length is $X_i(q)$.
Then for every $q$ we have
\[
    \big( \pmb{X}(q), \pmb{\mathfrak{N}}(q) \big) := \big( (X_1(q), X_2(q), \dots), (\mathfrak{N}_1(q), \mathfrak{N}_2(q), \dots) \big) \in \mathbb{U}_\searrow.
\]
We now state our second main result.
\begin{theorem}[Encoding the trajectories of the AMC]\label{thm:dynamicAMC2}
	The process
	\(
		\Big( \big( \pmb{X}(q), \pmb{\mathfrak{N}}(q) \big), q \ge 0 \Big) 
	\)
	has the law of the augmented multiplicative coalescent started from $(\pmb{x}, \pmb{0})$.
\end{theorem}
Theorem \ref{thm:dynamicAMC2} is proved in Section \ref{sec:sBFWencodesAMCdynamic}.

This construction and the existing results on the scaling limits of the simultaneous breadth-first walk (see \cite{multcoalnew}) are promising in the study of the eternal versions of the augmented multiplicative coalescent, as described in Section \ref{S:Scaling}.

The sBFW has already proved to be an efficient tool for studying the fluctuations of the size of the largest component in super-critical and barely super-critical Erd\H{o}s\,--\,Rényi\,--\,Stepanov random graph processes \cite{Corujo_Lemaire_Limic_2024}.

The framework discussed in this paper will be useful for a better understanding of the dynamic of the size of the connected components and the number of surplus edges in the super-critical Erd\H{o}s\,--\,Rényi\,--\,Stepanov random graph processes, as the next example illustrates.

\begin{example}[Erd\H{o}s--Rényi random graph]
This example can be seen as a extension of the Claim (Proposition 5, special case) in \cite[p.\ 2483]{multcoalnew}, as a corollary to Theorem \ref{thm:dynamicAMC2}.
Consider the case where $x_1 = x_2 = \dots = x_n = 1$ and $q = c/n$, with $c > 1$.
Denote $L_n(c)$ the size of the largest excursion of $B^{(1, \dots, 1), c/n}$ above zero, which occurs during the interval $[a_n(c), b_n(c)]$ and $\varLambda$ a homogeneous Poisson point process on $\mathbb{R}_+$.
Denote
\[
    S_n(c) = \varLambda \left( \frac{c}{n} \int_{a_n(c)}^{b_n(c)} B^{(1, \dots, 1), c/n} (s) \mathrm{d}s \right).
\]
 Theorem \ref{thm:dynamicAMC2} implies that the process 
\(
    \Big( \big( {L}_n(c), S_n(c) \big) , c \ge 0 \Big),
\)
follows the same law as the size and the number of surplus edges in the largest connected component of $(\mathcal{MG}^{(1/n, \dots, 1/n)} (c/n), c \ge 0 )$.
We are interested in computing the limits of
\[
    \frac{1}{n} L_n(c) \text{ and } \frac{1}{n} S_n(c),
\]
for every $c > 1$.
Notice that
\[
    Z^{(1, \dots, 1), c/n}(s) = \sum_{i = 1}^n \mathbb{1}_{\{\frac{\xi_i}{c/n} \le s\}} - s = n \cdot \tilde{Z}^{n,c} \left( \frac{s}{n} \right),
\]
where $(\xi_i)_{i = 1}^n$ are independent exponential rate $1$ random variables and
\[
    \tilde{Z}^{n,c}(s) = \frac{1}{n} \sum_{i = 1}^n \mathbb{1}_{\{\frac{\xi_i}{c} \le s\}} - s.
\]
By the Law of Large Numbers it is clear that
\[
    \tilde{Z}^{n,c}(s) \xrightarrow[n \to \infty]{\mathbb{P}} \Phi^c(s) = 1 - \mathrm{e}^{-c s} - s.
\]
Notice that $\Phi^c$ is concave and, for every $c > 1$, it has two zeros, one at $0$ and another at $\rho(c) > 0$. 
Hence, it is not difficult to derive from here the well-known asymptotic for the relative size of the giant component:
\[
    \frac{1}{n} L_n(c) \xrightarrow[n \to \infty]{\mathbb{P}} \rho(c).
\]
Let us denote by $\tilde{B}^{n,c}$ the reflected $\tilde{Z}^{n,c}$ above its past infimum.
It is easy to check that
\[
    S_n(c) = \Lambda \Big( c \int_{a_n(c)}^{b_n(c)} \tilde{B}^{n,c} \left( \frac{s}{n} \right) \mathrm{d}s \Big) = \Lambda \Big( c \, n  \int_{a_n(c)/n}^{b_n(c)/n} \tilde{B}^{n,c} (u) \mathrm{d}u \Big)
\]
Since $a_n(c)/n \to 0$ in probability, and $L_n(c) = b_n(c) - a_n(c)$, we get
\begin{equation}\label{eq:limit_surplus}
    \frac{1}{n} S_n(c) \xrightarrow[n \to \infty]{\mathbb{P}} c \int_0^{\rho(c)} \Phi^c(s) \mathrm{d}s = (c - 1) \rho(c) - c \frac{\rho(c)^2}{2}.
\end{equation}

We can easily check that the multi-graph $\mathcal{MG}^{(1, \dots, 1)}(c/n)$ has a $o(n)$ self-loops and multi-edges, for every $c \ge 0$.
Thus, the same asymptotics are true for the size and the number of surplus edges in the largest connected component of the super-critical Erd\H{o}s--Rényi random graph $G(n, c/n)$, with $c > 1$.
The relative asymptotic number of surplus edges \eqref{eq:limit_surplus} was obtained by Puhalskii \cite[Thm.\ 2.2]{Puhalskii2005}.
\hfill $\square$
\end{example} 
The forthcoming joint work with Nathanaël Enriquez, Gabriel Faraud and Sophie Lemaire uses the encoding of the AMC by the sBFW in the study of the joint fluctuations of the size of the largest component, the number of surplus edges and the number of connected components in super-critical Erd\H{o}s\,--\,Rényi random graphs.
Our Theorems \ref{thm:strictatqfull} and \ref{Theorem3.2B} also have a significant role in the forthcoming work \cite{limpla26}.

The coupling with surplus constructions exhibited in Sections \ref{S:take1} and \ref{S:take2} are intrinsic (up to randomization) to the simultaneous breadth-first walks.  
To the best of our knowledge, they also carry more detailed information than any of the previous surplus edge studies, see e.g.\ \cite{aldRGMC,bhamidietal2,bromar15,Dhara2017,Dhara2019}. 
In particular, provided that all the labels (positions) are kept for the surplus edges, the continuous-time random graph and the hereby
``enriched'' (simultaneous) breadth-first walks are equivalent, either in the sense of the marginal (see Lemma \ref{L:strictatq}) or the full distribution (see Theorem \ref{thm:strictatqfull}). 
Moreover, the coupling of Section \ref{S:take2} naturally motivates an extension to the multi-graph setting, which we link to \cite{bhamidietal2} in Section \ref{S:Multi-g}.
The dynamical encoding proposed in Section \ref{sec:sBFWencodesAMCdynamic} seems to be an optimal tool for studying the scaling limits of inhomogeneous random graphs in the near-critical regimes.
This is briefly discussed in Section \ref{S:Scaling}.
 
For general background on the random graph and the stochastic coalescence the reader is referred to \cite{aldous_survey, bertoin-fragcoal, bol-book, durrett-rgd,pitman-stflour}, and for specific as well as more recent references to \cite{multcoalnew}.
The edges in this paper will often be defined as oriented; however, when studying the global connectedness property in the resulting forest or (multi-)graphs, these orientations will not bear significance.

The rest of the paper is organized as follows.
In Section \ref{S:Coupl} we recall the breadth-first walk defined by the second author in \cite{multcoalnew} and we show how the surplus edges can be superimposed on this encoding.
A different encoding is discussed in Section \ref{S:take2}.
Although more complicated than the previous one, it produces a monotone forest representation of the multiplicative coalescent, and describes the relation between the surplus edges and the area above zero of the reflected simultaneous breadth-first walk.
In particular, we prove Theorem \ref{thm:encodingAMC}.
Our approach is based on ``ornamenting'' excursions (splitting the area under the excursion into ``interaction regions'') and careful tracking of surplus edges 
via Poisson marks in each region of interaction. 
In Section \ref{sec:ornamented} we provide a characterization of ornamented excursions. In Section \ref{sec:sBFWencodesAMCdynamic} we prove Theorem \ref{thm:dynamicAMC2} using a short argument, independent of the previous exposition in Section \ref{S:take2}.
Finally, in Section \ref{S:Scaling} we discuss some possible consequences of this paper for studying  eternal versions of the augmented multiplicative coalescent.

\section{The breadth-first walk and the surplus edges}
\label{S:Coupl}

\begin{figure}[h]
\centering
\includegraphics[scale=0.6]{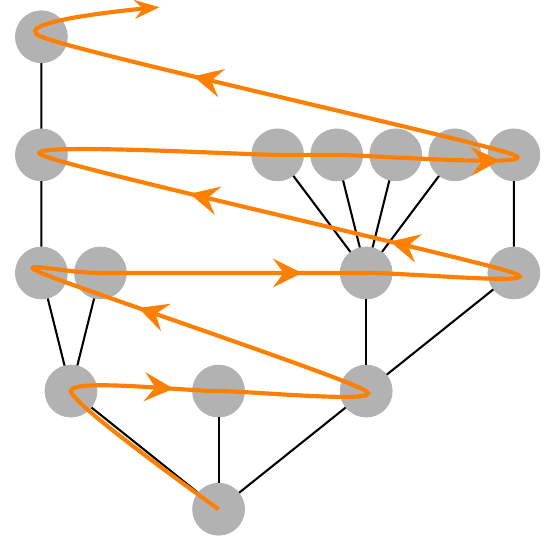}
\caption{Illustration of the breadth-first exploration of vertices in a finite rooted tree.
The vertices are explored according to the direction of the orange arrows.}
\label{F:one}
\end{figure}

Recall that ``breadth-first'' refers here to the order in which the vertices of a given connected
graph, or one of its spanning trees, are explored. 
Such exploration process starts at the root, visits all of its children, and these vertices become the first generation.
Then, it explores all the children of all the vertices from the first generation, and these vertices become the second generation. 
The exploration then turns to all the children of the second generation, and keeps going until all the vertices of all the generations are visited, or until forever if the tree is infinite.
Figure \ref{F:one} illustrates the breadth-first exploration of vertices in a finite rooted tree. 

In this section, random forests will conveniently span the components of a {\em coupled} {multiplicative coalescent}.  
In Section \ref{S:take1} these processes will be explored similarly to \cite{multcoalnew}, with a new feature: the non-spanning or {\em surplus} or {\em excess} edges will be recorded in addition. After that, in Section \ref{S:take2} another graph representation will be proposed in order to preserve the monotonicity of the graph-valued process.

\color{black}

\subsection{Breadth-first order induced forest and surplus}
\label{S:take1}
Recall the definition of $Z^{\pmb{x}, q}$ in \eqref{defZbxq}.
As recalled in the introduction, an important observation from \cite{multcoalnew} is that, simultaneously for all $q$, the  component masses  of the {multiplicative coalescent} started from $\pmb{x}$ and evaluated at time $q$ can be coupled to $Z^{\pmb{x},q}$ via a breadth-first walk construction. 
We next recall this construction.
See Section 2 in \cite{multcoalnew} for more details.

Given $a < b$ and an interval $[c, d]$ where $0 \le c < d$, define 
\(
(a, b] \oplus [c, d] := (a + c, b + d],
\) 
and denote by $|I|$ the Lebesgue measure of a Borel set $I$.
Typically $I$ will be an interval.
Let us also define $\xi^q_{(i)} := \xi_{(i)}/q$.

Algorithm \ref{alg:breadth-first} constructs the vector of component sizes associated to the simultaneous breadth-first walk defined by \eqref{defZbxq}.
Algorithm \ref{alg:breadth-first} takes as arguments the vector $\pmb{x}$ and the random vector  $\pmb{\xi}$ of exponentials, and constructs for a given $q > 0$, the edges $E$ of the random forest coupled to the simultaneous breadth-first walk, denoted by $\mathcal{F}_0(q)$.
In addition, the vector $\pmb{Y}(q)$ contains the sizes of all the trees (connected components) that were discovered up to the current step.

\renewcommand{\thealgocf}{A}
\begin{algorithm}[ht]

    $F_0^q \gets \{0\}$ \tcp*{initializing the load-free periods}

    $I_0^q \gets \{0\}$ \tcp*{initializing the listening periods}
	
 	$i \gets 1$ \tcp*{counter of the number of nodes}
	
	
 	$k \gets 0$ \tcp*{counter of the number of trees (connected components)}
	
 	$\pmb{Y}(q) \gets \mathbf{0}$ \tcp*{(unordered) vector of component sizes at time $q$}
	
 	\While{$i \leqslant \operatorname{len}(\pmb{x})$}{
		
 		$j \gets i$ \tcp*{auxiliary counter} \label{algo:aux_j}
		
 		$k \gets k + 1$ \label{algo:line_k}

         $F_k^q \gets \Big( \sup\big( I_{i-1}^q \big), \xi_{(i)}/q \Big]$ 
 		\tcp*{$k^{\text{th}}$ load-free period } \label{algo:def_Fk}

 		$I_i^q \gets ( \xi_{(i)}/q , \; \xi_{(i)}/q + x_{\pi_i} ]$ 
 		\tcp*{the listening period, $k^{\text{th}}$ leading vertex}
		
 		$Y_{k}(q) \gets Y_{k}(q) + x_{\pi_i}$
 		\tcp*{updating the component size}
		
 		\While{$ i \le j \le \mathrm{len}(\pmb{x})$}{ \label{algo:while_cycle}
			
 			\ForEach{ $\xi_{(r)}/q \in I_{i}^q \setminus I_{i-1}^q$ }{ \label{algo:for_cycle}
 				add edge $(\pi_r \to \pi_i)$ to $E$
\label{algo:add_parent}
                \tcp*{$\pi_r$ is a child of $\pi_i$}
				
 				$Y_{k}(q) \gets Y_{k}(q) + x_{\pi_r}$
 				\tcp*{updating the component size}
				
 				$j \gets j+1$
 			}
			
 			$I_{i+1}^q \gets I_i^q \oplus [0, x_{\pi_{i+1}}]$
 			\tcp*{updating the listening period}
			
 			$i \gets i + 1$ 	\tcp*{$j<i$ iff the connected component is exhausted}	
			
 		}
 	}
	
 	return $E$ and $\pmb{Y}(q)$

 	\caption{Implementation of the breadth-first walk algorithm, for given $\pmb{\xi}$ and $\pmb{x}$. 
    At each time $q \ge 0$ this constructs a random forest on $\{1,2, \dots, \operatorname{len}(\pmb{x})\}$ with edges $E$ and whose vector of component sizes $\pmb{Y}(q)$, ordered, evolves according to the multiplicative coalescent dynamic as $q$ varies. }  
 		\label{alg:breadth-first}

\end{algorithm}

Let us now explain this construction in words.
For $h = 1$ or $h$ equal to some of the values assigned when the Algorithm \ref{alg:breadth-first} visits line \ref{algo:line_k}, the corresponding $\pi_h$ is the root of a new tree in $\mathcal{F}_0(q)$ and for such values we have
\[
I_h^q = ( \xi_{(h)}/q , \xi_{(h)}/q + x_{\pi_h}],
\]
otherwise, when $\pi_h$ is not a root (here it must also be $h\leq \operatorname{len}(\pmb{x})$) we have 
\begin{equation}\label{def:Ihq}
    I_{h}^q = I_{h-1}^q \oplus [0, x_{\pi_{h}}].
\end{equation}
As a consequence,
$I^q_h \cap I^q_{h - 1} = \emptyset$ if and only if $\pi_h$ is a root of a new spanning tree.

The intervals $F^q_k$, for $k \ge 1$, defined in line \ref{algo:def_Fk} of Algorithm \ref{alg:breadth-first}, are the load-free periods.
Note that, for each $h \ge 1$, the dynamics ``listens'' for the children of $\pi_h$ during the interval $I^q_h \setminus I^q_{h - 1}$.
These intervals are called the {\em listening periods}.
The listening is implemented using the \emph{for loop} in line \ref{algo:for_cycle} of Algorithm \ref{alg:breadth-first}.
Furthermore, for each $i$ such that $\xi^q_{(i)} \in I^q_h \setminus I^q_{h - 1}$, we have that $\pi_h$ (resp.\ $\pi_i$) is the parent (resp.\ child) of $\pi_i$ (resp.\ $\pi_h$) in $\mathcal{F}_0(q)$.
In symbols, we write $\pi_i \to \pi_h$, thinking that any directed edge points to the parent, i.e.\ towards the root.

Finally, the above construction produces for each $q > 0$ a random forest $\mathcal{F}_0(q)$ with vertex set $[\operatorname{len}(\pmb{x})]$ and edge set $E$.
The output value of the variable $k$ is the number of trees in the forest $\mathcal{F}_0(q)$, and the output of Algorithm \ref{alg:breadth-first} is 
$\pmb{Y}(q) = (Y_1(q), Y_2(q), \dots, Y_k(q), 0, \dots)$, such that
\begin{equation}\label{eq:def_Mc_from_algorithm}
	\pmb{X}^{\pmb{x}}(q) := \operatorname{ord}\big((Y_1(q), Y_2(q), \dots, Y_k\textbf{}(q), 0, \dots) \big).
\end{equation} 
In words, it is a list of connected component masses of $\mathcal{F}_0(q)$, i.e.\ the masses of the trees obtained by adding up  the weights of member vertices (or equivalently, the masses of original particles given at time $0$). 
We naturally extend this definition at $q = 0$ to a trivial forest of $\operatorname{len}(\pmb{x})$ many single vertex trees ordered according to the order induced by $\pmb{\xi}$.

Note that Algorithm \ref{alg:breadth-first} simultaneously constructs the vector $\pmb{Y}(q)$ of component sizes at time $q$ and the random forest $\mathcal{F}_0(q)$.
This is done in line \ref{algo:add_parent} by assigning to each new vertex its parent.
Algorithm \ref{alg:breadth-first} is not intended to be optimal.
It is a choice, out of a number of possible pseudo-codes for the breadth-first walk algorithm.
We believe that it could improve the readers' understanding of the above construction.
Again, we refer the interested reader to Section 2 in \cite{multcoalnew} for more details on the construction of the random forest $\mathcal{F}_0(q)$ associated to the simultaneous breadth-first walk $Z^{\pmb{x}, q}$.

Define 
$$ B^{\pmb{x},q}(s):= Z^{\pmb{x},q}(s) - \inf_{u\leq s} Z^{\pmb{x},q}(u), \ s\geq 0, \ q >0.$$
Due to the above made observations, if and only if at time $s$ a new vertex is seen (or ``heard'') in $\mathcal{F}_0(q)$, the process $Z^{\pmb{x},q}$  (and therefore $B^{\pmb{x},q}$) makes an  upward jump of size equal to the mass of that vertex.
Let $k$ take a value assigned by the Algorithm \ref{alg:breadth-first} during a passage through line \ref{algo:line_k}.
The total sum of upward jumps of $B^{\pmb{x},q}$ during $T_k(q) := \operatorname{Cl}(I_{k}^q)$
is entirely compensated by the unit downward drift of $Z^{\pmb{x},q}$ (or $B^{\pmb{x},q}$) during $T_k(q)$. 
It is also easy to see that $B^{\pmb{x},q}(s)> 0$ in the interior of $T_k(q)$. 
Therefore, the $k^{\text{th}}$ excursion of $(B^{\pmb{x},q},\ q>0)$ above $0$ has the length precisely equal to the sum of the masses of all the blocks (vertices of the tree) explored during $T_k(q)$, i.e.\ $|T_k(q)| = Y_k(q)$.
It was initially noted in \cite{multcoalnew}, whose approach is rooted in \cite{aldRGMC,EBMC}, that for each fixed $q$, the ordered excursion lengths $\pmb{X}^{\pmb{x}}(q)$, as defined by \eqref{eq:def_Mc_from_algorithm}, have the {multiplicative coalescent} distribution, started from $\pmb{x}$ and evaluated at time $q$.
In particular, the exiting $K(q)\leq \operatorname{len}(\pmb{x})$ in the above algorithm, clearly equals the number of connected components at time $q$ of the coupled Aldous' continuous-time random graph. 
This type of result can be called \emph{static}, since it characterizes the marginal law of the process.

The advantage of the simultaneous breadth-first walk is that it also allows us to obtain a \emph{dynamical} result, characterizing the trajectories of the MC.
Proposition 5 in \cite{multcoalnew}, with foundations in 
\cite{Armendariz2001,Uribe2007},
shows the following stronger claim: $(\pmb{X}^{\pmb{x}}(q), q \ge 0)$ has the multiplicative coalescent law as a process in $q$, where the initial state at time $0$ is $\pmb{x}$.

One might wish to strengthen this in saying that the 
above coupling of $Z^{\pmb{x},q}$ and $\mathcal{F}_0(q)$ provides a bijective matching between the $k^{\text{th}}$ excursion of $B^{\pmb{x},q}$ above $0$, necessarily started at $\xi_{(h)}^q$ for some $h\in [\operatorname{len}(\pmb{x})]$, and a spanning (breadth-first explored) tree rooted at $\pi_h$ of the unique component of the continuous-time random graph which contains $\pi_h$ at time $q$. 
Some care is however needed here.

Indeed, while the mergers of different connected components or different subtrees of $\mathcal{F}_0(\cdot)$ arrive at precisely the multiplicative rate, the new edges arriving in $\mathcal{F}_0$ that correspond to those mergers always connect the root of one of the components (its excursion starts  at $\xi_{(h)}^q$ in $B^{\pmb{x},q -}$)  to the last visited or listed vertex in the previous component (its excursion is the one starting just before $\xi_{(h)}^q$ in $B^{\pmb{x},q  -}$).
In addition, within each connected component  the edges evolve according to a peculiar ``prune and reconnect'' rule, where vertices and subtrees are gradually ``moved closer'' to the root. 
Note that $\mathcal{F}_0(0)$ is the forest consisting in $\operatorname{len}(\pmb{x})$ isolated vertices, and for $q$ large enough, $\mathcal{F}_0(q)$ is the tree rooted at $\pi_1$ with $\operatorname{len}(\pmb{x})-1$ descendants.
In particular, the forest $(\mathcal{F}_0(q),\,q\geq 0)$ is not a monotone process with respect to the order induced by the subgraph relation.
Still, $\mathcal{F}_0(\cdot)$ is a random forest-valued process, whose tree masses evolve precisely according to \eqref{merge}, so it is a MC graph representation. 
We refer to any edge of $\mathcal{F}_0(q)$ as a {\em spanning edge at time} $q$. 

\subsection{Surplus on top of \texorpdfstring{$\mathcal{F}_0$}{F0}}\label{sec:surplus_edges}

There is a natural way to build the surplus or excess edges on top of $\mathcal{F}_0(q)$, for each $q$ separately, so that the distribution of the resulting graph is exactly the marginal of the continuous-time random graph at time $q$.
Since there are no loops or multi-edges in the continuous-time random graph $\big(\mathcal{G}^{\pmb{x}}(q), q \geq 0\big)$, we only need to check if there is a surplus edge between each pair of vertices which
are connected, but not by an edge in $\mathcal{F}_0(q)$.

\begin{figure}[h]
\begin{subfigure}[b]{0.3\textwidth}
	\centering
	\includegraphics[width=\textwidth]{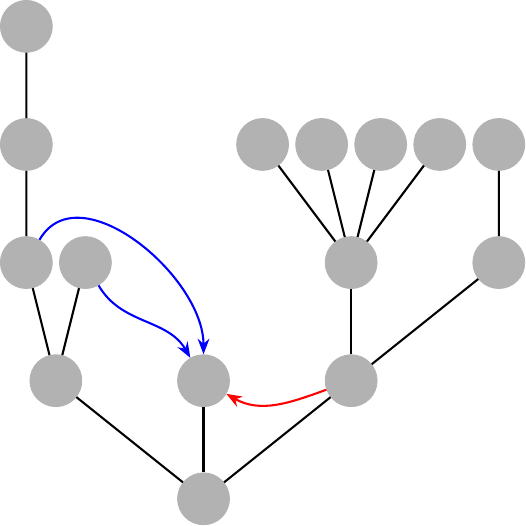}
	\caption{ }
	\label{F:two-a}
\end{subfigure}
\hfill
\begin{subfigure}[b]{0.3\textwidth}
	\centering
	\includegraphics[width=\textwidth]{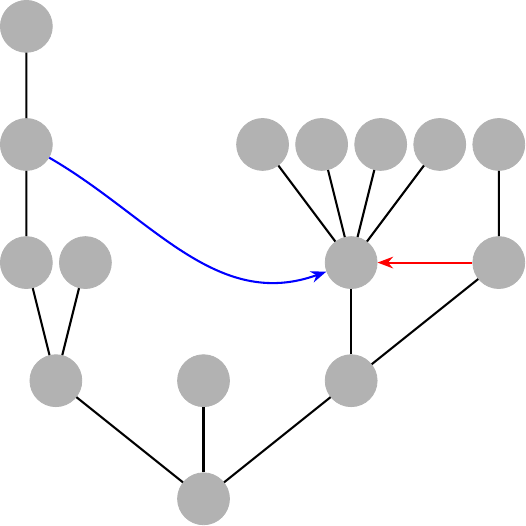}
	\caption{ }
	\label{F:two-b}
\end{subfigure}
\hfill
\begin{subfigure}[b]{0.3\textwidth}
	\centering
	\includegraphics[width=\textwidth]{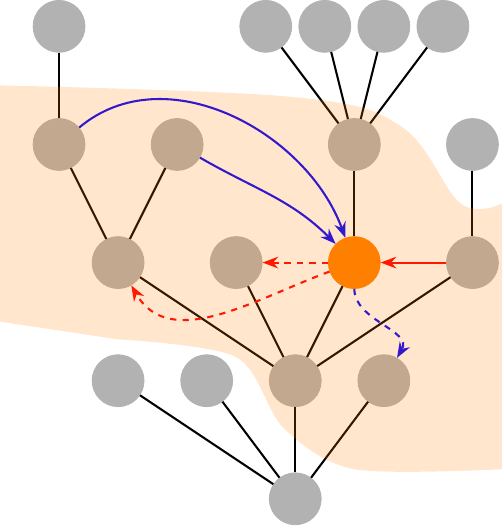}
	\caption{ }
	\label{F:three}
\end{subfigure}
\caption{
    Figures \subref{F:two-a} and \subref{F:two-b} show two different configurations and the surplus edges going to the third and the seventh vertex in the breadth-first order, respectively. 
    Possible surplus edges issued from vertices of the same generation, as described in \ref{item:surplus_same_gen} above, are indicated in red. Possible surplus edges issued from vertices of the next generation, as described in \ref{item:surplus_diff_gen}, are indicated in blue.
    Figure \subref{F:three} shows the influence region of the orange vertex, the eighth vertex in the breadth-first order.
    The red/blue coloring meaning is the same as in Figures \subref{F:two-a} and \subref{F:two-b}. 
    The arrows indicating a possible surplus edge between the orange vertex and a previously (resp.~subsequently) listed vertex are dashed (resp.~solid).
}
\label{F:two}
\end{figure}

Let us assume, without loss of generality, that we are looking for the surplus edges in the first tree in the random forest $\mathcal{F}_0(q)$.
Indeed, because of the lack of memory property of the exponential distribution, after exploring one or several components, the residual exponential times of the unexplored vertices have the same distribution as a family of ``new'' exponentials, independent from the past.
Suppose that for some $h\in [2,\operatorname{len}(\pmb{x})]$ we have that $\pi_h$ belongs to the first tree, that is $\xi_{(h)}^q \in I_{h-1}^q := (\xi_{(1)}^q, \xi_{(1)}^q + x_{\pi_1} + \dots + x_{\pi_{h-1}}]$.
Take $j \in [2,\operatorname{len}(\pmb{x})]$ such that $\xi_j^q \ge \xi_{(h)}^q$\footnote{Note that $\xi_j \ge \xi_{(h)}$ can be equivalently written as $\pi^{-1}_j \ge h$, where $\pi^{-1}$ is the inverse permutation of $\pi$, in the sense that $\pi_i = k$ if and only if $\xi_{(i)} = \xi_{k}$.}, i.e.\ $j$ comes after $\pi_{h}$ according to the breadth first order.
Then, on
\[
	\{\xi_j^q \in \big(\xi_{(1)}^q + x_{\pi_{1}} + \dots + x_{\pi_{h-1}}, \xi_{(1)}^q + x_{\pi_{1}} + \dots + x_{\pi_{h-1}} + x_{\pi_{h}}\big]	\}
\]
the edge $j \to \pi_h$ belongs to $\mathcal{F}_0(q)$, and the above event occurs with probability $1 - \mathrm{e}^{-q \, x_{\pi_h} \, x_j}$.
This is a consequence of the breadth-first search for descendants  (see \cite[\S\ 2]{multcoalnew}) and of the lack of memory of the exponential distribution.
Indeed, let us denote
$$
	\mathcal{Z}^q_s = \sigma\left( \{ (\xi_i^q > u): i \in [\operatorname{len}(\pmb{x})] \}, u \le s \right).
$$
Hence, $(\mathcal{Z}_s^q)_{s \ge 0}$ is the filtration generated by the arrivals of $(\xi_1^q, \xi_2^q, \dots)$.
Recall that
\[
    \{\xi_j^q > \sup(I_{h-1}^q)\} = \{\xi_j^q > \xi_{(1)}^q + x_{\pi_1} + \dots + x_{\pi_{h-1}}\} \in \mathcal{Z}_{\sup(I_{h-1}^q)}^q.
\]
Then, we have
\begin{equation*} 
	\mathbb{P}\left(\xi_j^q > \sup(I_{h}^q) \mid \mathcal{Z}_{\sup(I_{h-1}^q)}^q \right) \mathbb{1}_{\left\{ \xi_j^q > \sup(I_{h-1}^q) \right\} } = \mathrm{e}^{- q \, x_{j} \, x_{\pi_h}} \mathbb{1}_{\left\{ \xi_j^q > \sup(I_{h-1}^q) \right\} } \;\; \text{ a.s.}
\end{equation*}
Note that $\mathrm{e}^{- q \, x_j \, x_{\pi_h}}$ is precisely the probability that vertices $j$ and $\pi_{h}$ are not connected at time $q$ in the Aldous' continuous-time random graph dynamic.
So the above breadth-first walk encoding procedure does not only yield the spanning edges according to the random graph dynamics, but also prohibits existence of the edges $j \to \pi_h$, where $j$ is any vertex listed after the children of $h$ in the breadth-first order induced by the algorithm.
Note that, if $\xi_{j}^q \in ( \xi_{(h)}^q , \sup (I_{h-1}^q) ]$, then $j$ is necessarily another vertex in the same tree as $\pi_h$, encountered after $\pi_h$ and before any child of $\pi_h$.
Thus, due to the above breadth-first exploration order, see Figure \ref{F:one}, there are only two types of configurations where a surplus edge $j\to \pi_h$ is possible:
\begin{enumerate}[label=(\roman*)]
\item $j$ is listed after $\pi_h$, in the same generation as $\pi_h$, or
 \label{item:surplus_same_gen}
\item $j$ is from the next generation (then necessarily a child of some vertex $\pi_k$ with $k<h$, where $\pi_k$ is from the same generation as $\pi_h$). \label{item:surplus_diff_gen}
\end{enumerate}
Similar observations in the setting of discrete breadth-first walks appear in Aldous \cite{aldRGMC} and Spencer \cite{Spencer97}.
See Figures \ref{F:two-a} and \ref{F:two-b} for an illustration of these two cases. 

The oriented surplus edges with arrows pointing to $\pi_h$ are possible only from those indices $j$ satisfying \ref{item:surplus_same_gen} or \ref{item:surplus_diff_gen}.
Indeed, if $j$ is encountered by the BF encoding after all the children of $\pi_{h}$, then there is no spanning or surplus edge connecting $j$ and $\pi_{h}$. 
If $\xi_j^q$ enters $I_h^q\setminus I_{h-1}^q$, as noted above, this means that  $j$ is a child of $\pi_{h}$ (which prohibits existence of a surplus edge connecting them). 
Finally, if $\xi_j^q \in I_{h-1}^q \setminus [0,\xi_{(h)}^q]= ( \xi_{(h)}^q , \sup (I_{h-1}^q) ]$ then either \ref{item:surplus_same_gen} or \ref{item:surplus_diff_gen} are true, and a surplus edge connecting $j$ and $\xi_{(h)}$ is not excluded.
So, we need an external source of randomness, independent from the algorithm, in order to determine whether such surplus edges exist or not.

Figure \ref{F:three} shows the ``surplus influence region''  for a typical non-root vertex.
Note that the surplus influence region of a given vertex $v$ consists of all vertices attached to the tree (strictly) between the parent of $v$ and the first child of $v$, in the sense of the order
induced by the breadth-first algorithm.

The presence or absence of surplus edges should be compatible with the random graph dynamic.
For each $l,h \in [\operatorname{len}(\pmb{x})]$, let $\zeta^{l,h}$ be a Poisson process of marks arriving at the rate $x_{l}$.
The processes $(\zeta^{l, h})_{l, h}$ are independent over $l$ and $h$, and are independent of $\pmb{\xi}$. 
Consider those (and only those) $\zeta^{l,h}$ such that $\xi_{l}^q \in I_{h-1}^q\setminus [0, \xi_{(h)}^q]$.     
For any such $l$ draw the (red or blue, see Figure \ref{F:two}) surplus edge connecting $l$ and $\pi_h$ if and only if 
\begin{equation}
    \label{eq:jump_in_zeta}
    \zeta^{h,l}[0, q \cdot x_{\pi_h}] \geq 1.
\end{equation} 
It is clear that this edge appears with probability $1-\mathrm{e}^{-q \, x_{l} \, x_{\pi_h}}$, independently of everything else. 
Any other edge at time $q$ is also present with probability $1-\mathrm{e}^{-q \, m_1 \, m_2}$, where $m_1$ and $m_2$ are the masses of the two vertices (see above discussion and Section 2 in \cite{multcoalnew} for more details).
 
Denote by $G^{\pmb{x}}(q)=(V,E(q))$ the resulting random graph, where $E(q)$ is the union of the spanning and the surplus edges at time $q$.
Recall the Aldous continuous-time random graph $\mathcal{G}^{\pmb{x}}(q)$. 
One can record the just made observations as follows.
\begin{lemma}
\label{L:strictatq}
For each $q \geq 0$ and every finite vector $\pmb{x}$, the random graph $G^{\pmb{x}}(q)$ has the same distribution as $\mathcal{G}^{\pmb{x}}(q)$.
\end{lemma}

\begin{remark}[Joint intensity of the Poisson processes] \label{rmk:joint_intensity}
Suppose that one is only interested in counting the surplus edges in various connected components of the random graph, without keeping track of their exact location (meaning the vertices which it connects).
Note that
the joint (total) intensity of all the Poisson marking processes $\zeta^{l,h}$ for a fixed value of $h$ is 
\begin{align}
\sum_{l\,:\,\xi_l^q \in I_{h-1}^q\setminus [0, \xi_{(h)}^q]} x_l &= \left(\sum_{s\in  I_{h-1}^q} \Delta Z^{\pmb{x}, q}(s) - |I_{h}^q| \right) = \Big(B^{\pmb{x},q}\big(\sup(I_{h-1}^q)\big)- x_{\pi_h}\Big), \label{Etotalinten}
\end{align}
where $\Delta Z^{\pmb{x}, q}(s) := Z^{\pmb{x}, q}(s) - Z^{\pmb{x}, q}(s-)$, and where we recall that $B^{\pmb{x},q}$ equals $Z^{\pmb{x},q}$ reflected above past minima. 
Indeed, the first identity in \eqref{Etotalinten} is a consequence of the following identities
\[
\sum_{l\,:\,\xi_{l}^q \in I_{h-1}^q} x_{l} =\sum_{s\in  I_{h-1}^q} \Delta Z^{\pmb{x}, q}(s)  \;\; \text{ and } \;\; 
\sum_{l \,:\, \xi_l^q \in [\xi_{(k)}^q, \xi_{(h)}^q]} x_l = |I_{h}^q| = x_{\pi_{k}} + \dots + x_{\pi_{h}},
\]
where $\pi_k$ is the root of the tree including $\pi_h$ as one of its vertices.
To check the second identity in \eqref{Etotalinten} note that, due to the drift in the definition of $Z^{\pmb{x}, q}$, as well as the definition of the reflection map transforming $Z^{\pmb{x}, q}$ into $B^{\pmb{x}, q}$, we have
\[
	B^{\pmb{x},q}\Big(\sup(I_{h-1}^q)\Big) = \sum_{s \in I_{h-1}^q} \Delta Z^{\pmb{x}, q}(s) - |I_{h-1}^q|.
\]
Recall that typically $Z^{\pmb{x},q}$  does not  attain its global infimum on $[0,\sup(I_h^q)]$. 
However, its infimum on the interval $[0,\sup(I_h^q)]$ is always equal to the left limit  $Z^{\pmb{x},q}(\xi_{(k)}^q-)$ of $Z^{\pmb{x},q}$ evaluated at the left endpoint $\inf(I_h^q)=\xi_{(k)}^q$ of $I_h^q$. 
Recall that $\pi_k$ is the root of the corresponding tree (as above).

 When we multiply \eqref{Etotalinten} by the elapsed time $q \cdot x_{\pi_h}$
 (see \eqref{eq:jump_in_zeta} above)
we obtain the mean $A_h$
of the Poission random variable which specifies the number of \emph{unlimited surplus edges} (meaning that multiple edges are allowed)   incident to  $\pi_h$:
\[
    A_h = q \cdot x_{\pi_h} \cdot \sum_{l>h,\, \xi_{(l)}^q\in I_{h-1}^q} x_{\pi_l} = q \Big( x_{\pi_h} \cdot B^{\pmb{x},q}\big(\sup(I_{h-1}^q)\big)- x_{\pi_h}^2\Big).
\]
Hence, the total number of unlimited surplus edges within a connected component $C$ is a Poisson random variable with mean
\[
    \sum_{h: \, \pi_h \in C} A_h = \sum_{h: \, \pi_h \in C} q \Big( x_{\pi_h} \cdot B^{\pmb{x},q}\big(\sup(I_{h-1}^q)\big)- \big( x_{\pi_h} \big)^2 \Big).
\]

\end{remark}

The next remark is a digression which could help an interested reader in understanding the consequence of the computations in Remark \ref{rmk:joint_intensity} on the surplus edge analysis in the near-critical regime.
Our second construction of surplus edges (see Section \ref{S:take2})
leads to a precise formulation and analysis of an analogous property.
In Section \ref{S:Scaling} we discuss some interesting consequences in the near-critical regime.

\begin{remark}[The link between the surplus edges and the area below the curve] \label{rmk:joint_intensity2}
Suppose that $C_i^{(n)}:[a_i^{(n)},b_i^{(n)}] \mapsto \mathbb{R}_+$ is the $i^{\text{th}}$ longest excursion of $B^{\pmb{x}^{(n)},q_n(t)}$ above $0$,
where $\pmb{x}^{(n)}$ satisfies conditions (5.1)-(5.3) in \cite{multcoalnew} (which are also given as \eqref{eq:hypo1}, \eqref{eq:hypo2}, \eqref{eq:hypo3} below) and that $q_n(t) = t + 1/\sigma_2^{(n)}$, where $\sigma_2^{(n)} = \sum (x_i^{(n)})^2$.

    In the scaling limit (see \cite[Prop.\ 1.6]{multcoalnew} or Section \ref{S:Scaling} below) it is natural to expect
    \begin{equation*}
        \label{eq:natural_expectation}
        \lim\limits_{n \to \infty} q_n(t) \sum_{h: \, \pi_h \in C_i^{(n)}} \Big( x_{\pi_h}^{(n)} \cdot B^{\pmb{x}^{(n)},q_n(t)}\big(\sup(I_{h-1}^q)\big) -  \big( x_{\pi_h}^{(n)} \big)^2\Big) = \int_{a_{i-1}}^{a_i} B^{\kappa, t, \pmb{c}}(s) \mathrm{d}s,
    \end{equation*}
    where $a_{i-1}$ and $a_i$ are respectively the left and the right endpoints of the $i^{\text{th}}$ longest excursion of $B^{\kappa, t, \pmb{c}}$ (defined in \eqref{eq:def_B^kappa^tau^c}) above zero.
    Indeed, 
    \[
        \sum_{h: \, \pi_h \in C_i^{(n)}} x_{\pi_h}^{(n)} \cdot q_n(t) B^{\pmb{x}^{(n)},q_n(t)}\big(\sup(I_{h-1}^q)\big)
    \]
    is a Riemann sum of $q_n(t) B^{\pmb{x}^{(n)},q_n(t)}$ restricted to the domain $[a_i^{(n)}, b_i^{(n)}]$ of $C_i^{(n)}$.
    In the scaling limit $ q_n(t) B^{\pmb{x}^{(n)},q_n(t)}$  converges towards $B^{\kappa, t, \pmb{c}}$ in the sense of the $J_1$ topology \cite[Prop.\ 6]{multcoalnew}.
    In addition, the error term
	\[
       q_n(t) \sum_{h: \, \pi_h \in C_i^{(n)}} (x_{\pi_h}^{(n)})^2 = 
       t \sum_{h: \, \pi_h \in C_i^{(n)}} (x_{\pi_h}^{(n)})^2 
       + \frac{1}{\sigma_2^{(n)}} \sum_{h: \, \pi_h \in C_i^{(n)}} (x_{\pi_h}^{(n)})^2
    \]
    is negligible  (again due to assumptions \eqref{eq:hypo1}--\eqref{eq:hypo3}).
Indeed, the first term on the RHS goes to zero when $n$ is large, because of \eqref{eq:hypo3}.
The second term can be controlled using the Cauchy\,--\,Schwarz inequality
\begin{align*}
    \frac{1}{\sigma_2^{(n)}} \sum_{h: \, \pi_h \in C_i^{(n)}} (x_{\pi_h}^{(n)})^2 &= \sum_{h: \, \pi_h \in C_i^{(n)}} \frac{(x_{\pi_h}^{(n)})^{3/2}} {(\sigma_2^{(n)})^{3/2}} \times \big( \sigma_2^{(n)}  \times x_{\pi_h}^{(n)}  \big)^{1/2} \\ 
    &\le \left( \sum_{i \ge 1}  \frac{(x_{i}^{(n)})^{3}} {(\sigma_2^{(n)})^{3}} \right)^{1/2} \Big( \sigma_2^{(n)} \times X_i^{(n)} \Big)^{1/2}.
\end{align*}
Using \eqref{eq:hypo1} we get that the first factor in the previous display converges to $(\kappa + \sum_i c_i^3 )^{1/2}$, which is a 
 finite quantity.
In addition, the second factor converges to zero in law, again by \eqref{eq:hypo3} and the fact that  $X_i^{(n)}=b_i^{(n)}- a_i^{(n)}$ converges in law to to the length of the $i^{\text{th}}$ longest excursion of $B^{\kappa, t, \pmb{c}}$.
    
In the setting of simple graphs, where there is at most one edge between any pair of vertices (see Algorithm \ref{alg:breadth-first} and \eqref{eq:jump_in_zeta}), the above reasoning shows that the number of surplus edges in the component which corresponds to $C_i^{(n)}$ is a sum of independent Bernoulli random variables with parameters converging to zero as $n \to \infty$.
Moreover the sum of the parameters over all the vertices in this component converges to $\int_{a_{i-1}}^{a_i} B^{\kappa, t, \pmb{c}}(s) \mathrm{d}s$.
Hence, using Le Cam's approximation of the sum of independent Bernoulli random variables by a Poisson distribution \cite[Thm.\ 2]{LeCam1960} we get that the number of limited surplus edges converges also to a Poisson random variable with parameter $\int_{a_{i-1}}^{a_i} B^{\kappa, t, \pmb{c}}(s) \mathrm{d}s$.

An analogous result   was derived in \cite[\S \, 2.2]{aldRGMC} in the  setting of critical Erd\H{o}s\,--\,R\'enyi graphs and the Aldous' standard multiplicative coalescent. 
\end{remark}

\section{Monotone forest representation and surplus edges}
\label{S:take2}

We next revise and enrich (via an additional randomization) the coupling algorithm just described. 
Instead of the random forest process $(\mathcal{F}_0(q),\, q>0)$, another forest-valued process $(\mathcal{F}_1(q),\, q>0)$ is described, so that it is motonone ($\mathcal{F}_1(q_1)$ is a subgraph of $\mathcal{F}_1(q_2)$ when $q_1 \le q_2$).
The process $\mathcal{F}_1$ is linked to the evolution of the Aldous' random graph in a sense that we now explain.

As before, let us fix a finite vector $\pmb{x}$.
Let $\big( \mathcal{G}^{\pmb{x}}(q), q \ge 0 \big)$ be the inhomogeneous random graph started from $\operatorname{len}(\pmb{x})$ many isolated nodes of masses given by $\pmb{x}$.
Let $\big( \pmb{X}^{\pmb{x}}(q), q \ge 0 \big)$ be the MC coupled to $\big( \mathcal{G}^{\pmb{x}}(q), q \ge 0 \big)$, i.e.\ the process 
recording the size of the connected components in $\mathcal{G}^{\pmb{x}}(q)$, for every $q \ge 0$.
Consider now a forest-valued process $\big( \mathcal{T}^{\pmb{x}}(q), q \ge 0 \big)$ such that,  $\big( \mathcal{T}^{\pmb{x}}(q), q \ge 0 \big)$ only jumps when $\big( \pmb{X}^{\pmb{x}}(q), q \ge 0 \big)$ does, i.e.\ when two connected components of $\mathcal{G}^{\pmb{x}}(q-)$ merge into one new component of $\mathcal{G}^{\pmb{x}}(q)$.
A transition step of $\big( \mathcal{T}^{\pmb{x}}(q), q \ge 0 \big)$ consists of appending one edge to the graph $\mathcal{T}^{\pmb{x}}(q-) $.
This new edge in $\mathcal{T}^{\pmb{x}}(q)$ is precisely the new edge appearing in $\big( \mathcal{G}^{\pmb{x}}(q), q \ge 0 \big)$ at the moment of the corresponding jump, i.e. it is the edge which produces the merger.
One can deduce from the previous construction that $\big( \mathcal{T}^{\pmb{x}}(q), q \ge 0 \big)$ satisfies the following properties:
\begin{enumerate}[label=\textbf{M.\arabic*}]
    \item $\big( \mathcal{T}^{\pmb{x}}(q), q \ge 0 \big)$ is a forest-valued Markov process, \label{property:M1}
    \item $\big( \mathcal{T}^{\pmb{x}}(q), q \ge 0 \big)$ is monotone in the sense of the inclusion of sets, \label{property:M2}
    \item for every $q \ge 0$, $\mathcal{T}^{\pmb{x}}(q)$ is a spanning forest of $\mathcal{G}^{\pmb{x}}(q)$, meaning that each tree in $\mathcal{T}^{\pmb{x}}(q)$ is a spanning tree of some connected component of $\mathcal{G}^{\pmb{x}}(q)$. \label{property:M3}
    \item given the past and present configurations $\sigma (\mathcal{T}^{\pmb{x}}(z): z\leq q )$, and knowing that the random graph at time $q$ (resp.~the random forest $\mathcal{T}^{\pmb{x}}(q)$) has $N$ connected components (resp.~spanning trees), which induce a partition 
    \[
    \{i_1,i_2,\ldots,i_{k_1}\}, \{i_{k_1 + 1},i_{k_1 + 2},\ldots,i_{k_2}\}, \dots, \{i_{k_{N-1}+1},i_{k_{N-1}+2},\ldots,i_{\operatorname{len}(\pmb{x})}\}
    \]
    on $[\operatorname{len}(\pmb{x})]$ (without knowing the breadth-first order induced by $\pmb{\xi}$, but knowing that the mass of $i$ equals $x_i$), the edge between any two $h$ and $l$, which belong to different trees of $\mathcal{T}^{\pmb{x}}(q)$, arrives at rate $x_h \cdot x_l$.
\label{property:M4a}
\end{enumerate}
Several random forest-valued processes can be constructed such that properties \ref{property:M1}, \ref{property:M2} and \ref{property:M3} are satisfied.
For example, we can construct a trivial such process in the following way: at the moment when its connected component merges with another (previously listed) connected component,
the leading vertex
$\pi_h$ (it is the vertex which serves as the root of its component just prior to the merger time) always gets connected by an edge to $\pi_{h-1}$.
This unsophisticated procedure constructs a monotone forest-valued process that evolves starting from the forest of $\operatorname{len}(\pmb{x})$ many trivial trees (i.e.\ $\mathcal{F}_0(0)$), and ends with a tree of $\operatorname{len}(\pmb{x})$ many generations, with a single vertex in each generation.
However, the edge arrivals in this construction do not satisfy \ref{property:M4a}. 
Indeed, if two connected components, composed of $\{i_1,i_2,i_3\}$ and $\{i_4,i_5\}$, merge at a certain time $q$ via an edge connecting $i_2$ and $i_5$, then (given all the available information) we know that starting from time $q$ both $i_2$ and $i_5$ are interior vertices within the tree containing $\{i_1,i_2,i_3,i_4,i_5\}$, and therefore it is not possible for either of these vertices to connect via an edge to any other vertex at a future merging event. 

The process $\big( \mathcal{T}^{\pmb{x}}(q), q \ge 0 \big)$ can be thought of as a \emph{dynamical spanning forest} coupled to the random graph $\big( \mathcal{G}^{\pmb{x}}(q), q \ge 0 \big)$. 
Our goal in this section is to construct a forest valued process $\big( \mathcal{F}_1(q), q \ge 0 \big)$ coupled to the simultaneous breadth-first walks $(Z^{{\pmb{x}}, q}, q \ge 0)$, such that it has the same distribution as $\big( \mathcal{T}^{\pmb{x}}(q), q \ge 0 \big)$.

Let us take a closer look at the dynamic of the process $\big( \mathcal{T}^{\pmb{x}}(q), q \ge 0 \big)$.
This Markov process jumps according to the multiplicative coalescent rates $X_i \cdot X_j$, where $X_i$ and $X_j$ are the total masses of two connected components (trees) in $\mathcal{T}^{\pmb{x}}(q)$.
Conditionally on the fact that the components $C_i = \{i_1, \dots, i_k\}$ and $C_j = \{j_1, \dots, j_m\}$ are merging, where $|C_i| = X_i = x_{i_1} + \dots + x_{i_k}$ and $|C_j| = X_j = x_{j_1} + \dots + x_{j_m}$, and that the merging time is $Q$,
the new edge in $ \mathcal{T}^{\pmb{x}}(Q)$ is $\{i_s, j_t\}$ with probability
\[
    \frac{x_{i_s} }{ X_j} \frac{x_{j_t}}{X_j}, \text{ for every } i_s \in C_i \text{ and } j_t \in C_j.
\]
Indeed, suppose that, for 
$ i_s \in C_i$, $j_t\in C_j$, $\zeta^{i_s, j_t}$ is an exponential random variable with rate $ x_{i_s} \cdot x_{j_t}$, independently over different $s$ and $t$.
Then, according to the random graph dynamic and elementary properties of the exponential distribution, on the event $\{Q=q\}$ the edge $\pi_{i_s} \to \pi_{j_t}$ is the one created by the merger of $C_i$ and $C_j$  with probability
\begin{align*}
	\mathbb{P}\left[ \zeta^{i_s,j_t} = \min_{k \in C_i, l \in C_j} \zeta^{k,l} \,  \Big| \, \min_{k \in C_i, l \in C_j} \zeta^{k,l}, \mathcal{T}^{\pmb{x}}(q-) \right] &\equiv \frac{q \cdot x_{i_s} x_{j_t} }{ \sum\limits_{k \in C_i, l \in C_j} q \cdot x_{k} x_{l} } = \frac{x_{i_s}}{ X_i} \frac{ x_{j_t} }{ X_j}.
\end{align*}
This is equivalent to picking one vertex in each component in a size-biased way, and connecting them by an edge.
The arrival of this edge reduces the total number of trees by one
(it creates one larger tree from two smaller trees).

We now describe the construction of $(\mathcal{F}_1(q), q \ge 0)$.
The initial forest $\mathcal{F}_1(0)$ is again trivial, and therefore equal to $\mathcal{F}_0(0)$.
It is the set of vertices $[\operatorname{len}(\pmb{x})]$, without any edge.
During a strictly positive (random) interval of time $\mathcal{F}_1(\cdot)$ will remain $\mathcal{F}_1(0)$.
Recall the definition of the listening periods \eqref{def:Ihq}.
At some (stopping) time $Q_1>0$ the first connection is established between $\pi_{L}$ and $\pi_{L-1}$, where $L$ is the first value in $[\operatorname{len}(\pmb{x})]$ such that $\xi_{(L)}^{Q_1}\in I_{L - 1}^{Q_1} = (\xi^{Q_1}_{(L - 1)}, \xi^{Q_1}_{(L - 1)} + x_{\pi_{L-1}}]$.
In particular, $\xi_{(h)}^{Q_1} \not \in I_{h-1}^{Q_1}$ for any $h\in [2,\operatorname{len}(\pmb{x})] \setminus \{L\}$, almost surely.
At time $Q_1$ both $\mathcal{F}_0$ and $\mathcal{F}_1$ make the same jump: the new edge $\pi_{L}\rightarrow \pi_{L-1}$ appears. 
After that, during an interval of time of positive random length, $\mathcal{F}_1$ stays equal to $\mathcal{F}_1(Q_1)$, and eventually a new connection occurs at some time $Q_2>Q_1$. 
The difference between $\mathcal{F}_1$ and $\mathcal{F}_0$ may be visible already at time $Q_2$. 
Indeed, for the sake of illustration suppose that either $\xi_{(L+1)}^{Q_2}$ enters $I_{L}^{Q_2} = ( \xi^{Q_2}_{(L-1)}, \xi^{Q_2}_{(L-1)} + x_{\pi_{L-1}} + x_{\pi_{L}} ]$, or $\xi_{(L-1)}^{Q_2}$ enters $I_{L-2}^{Q_2} = ( \xi^{Q_2}_{(L-2)}, \xi^{Q_2}_{(L-2)} + x_{\pi_{L-2}} ]$ at this very moment. 
As already noted, in  $\mathcal{F}_0(Q_2)$ we will add the edge $\pi_{L + 1}  \rightarrow \pi_{L}$ in the former case, and $\pi_{L-1}  \rightarrow \pi_{L-2}$ in the latter case.
Also note that in $\mathcal{F}_0$ both $\pi_{L}$ and $\pi_{L - 1}$ will eventually become children of the same vertex $\pi_{h}$ for some $h \in \{1,\dots, L -2\}$ and this $h$ will change until finally becoming equal to $\pi_1$.
However, given $\{ \xi_{(L+1)}^{Q_2} \in I_{L}^{Q_2}  \}$ the new edge in $\mathcal{F}_1(Q_2)$ is chosen as
\begin{itemize}
	\item $\pi_{L+1}  \rightarrow \pi_{L}$, with probability $\displaystyle \frac{x_{\pi_{L}}}{x_{\pi_{L-1}}+ x_{\pi_{L}}}$, or
	\item $\pi_{L + 1}  \rightarrow \pi_{L - 1}$, with the remaining probability $\displaystyle \frac{x_{\pi_{L - 1}}}{x_{\pi_{L - 1}} + x_{\pi_{L}}}$.
\end{itemize}
Similarly, given  $\{ \xi_{(L-1)}^{Q_2} \in I_{L-2}^{Q_2}  \}$, the new edge in $\mathcal{F}_1(Q_2)$ is chosen as
\begin{itemize}
	\item $\pi_{L - 1}  \rightarrow \pi_{L - 2}$, with probability $\displaystyle \frac{x_{\pi_{L - 1}}}{x_{\pi_{L - 1}}+ x_{\pi_{L}}}$, or
	\item $\pi_{L }  \rightarrow \pi_{L - 2}$, with the remaining probability $\displaystyle \frac{x_{\pi_{L}}}{x_{\pi_{L-1}} + x_{\pi_{L}}}$.
\end{itemize}
Once it appears in $\mathcal{F}_1$, this and any other edge will remain in $\mathcal{F}_1$ forever after.
It is important to note that these transitions  
are compatible with the evolution of $(\mathcal{T}^{\pmb{x}}(q), q \ge 0)$.

The complete construction of $\mathcal{F}_1$ is as follows: $\mathcal{F}_1(0)=\mathcal{F}_0(0)$; for $q>0$, $\mathcal{F}_1(q)=\mathcal{F}_1(q-)$ unless $q$ is such that the number of components (trees) in $\mathcal{F}_0$ at time $q$ decreases by $1$.
The latter happens if and only if for some $j \in [\operatorname{len}(\pmb{x})]$, $\xi_{(j)}^{q-} \not \in I_{j-1}^{q-}$ (recall that $\pi_j$ is the root of a spanning tree in $\mathcal{F}_0(q-)$, or equivalently, in $\mathcal{F}_1(q-)$) and $\xi_{(j)}^{q} \in I_{j-1}^{q}$.
Assume that the two components that merge at time $q$ in $\mathcal{F}_0(q)$ are
\(
	\{ \pi_i, \pi_{i+1}, \dots, \pi_{j-1} \} \text{ and } \{ \pi_j, \pi_{j+1}, \dots, \pi_{j + k} \}.
\)
For such $q$, let $\mathcal{F}_1(q)$ inherit all the edges of $\mathcal{F}_1(q-)$, and in addition 
\begin{equation}
\label{Enewedge}
	\text{draw a new edge } {L_1(q)} \rightarrow {L_2(q)} \mbox{ in } \mathcal{F}_1(q),
\end{equation}
where $L_1(q)$ (resp.~$L_2(q)$) is chosen at random (in the size-biased way) from the vertices of  $\{ \pi_j, \pi_{j+1}, \dots, \pi_{j + k} \}$ (resp.~$\{ \pi_i, \pi_{i+1}, \dots, \pi_{j-1} \}$). 
More precisely, conditionally on $\{\mathcal{F}_1(s),\, s<q\}$ and $\{\mathcal{F}_0(s),\, s\leq q\}$, let $L_1(q)$ equal $\pi_s$ with probability $x_{\pi_s}/(x_{\pi_j}+ x_{\pi_{j+1}} + \dots + x_{\pi_{j+k}} )$ for each $s \in \{j, j+1, \dots, j+k\}$, and $L_2(q)$ equal $\pi_l$ with probability $x_{\pi_l}/(x_{\pi_i}+ x_{\pi_{i+1}} + \dots + x_{\pi_{j-1}} )$ for each $l \in \{i, i+1, \dots, j-1\}$, where these choices are conditionally independent.
Finally assign the new edge via \eqref{Enewedge}.

It is clear from the just presented construction that $\mathcal{F}_1$ is a monotone forest-valued process, and also that, almost surely, for each $q$ the trees in $\mathcal{F}_1(q)$ are composed of exactly the same vertices as the (corresponding) trees in $\mathcal{F}_0(q)$ (or equivalently, $i$ and $j$ are connected in $\mathcal{F}_1(q)$ if and only if they are connected in $\mathcal{F}_0(q)$).  
Moreover, we obtain the following.
\begin{lemma}
    $(\mathcal{F}_1(q), q \ge 0)$ follows the same distribution as the dynamical spanning forest $(\mathcal{T}^{\pmb{x}} (q), q \ge 0)$.
\end{lemma}

\begin{remark}
The just described way of attaching edges in $\mathcal{F}_1$ is the most natural (if not the only) choice for a monotone MC forest representation coupled with the breadth-first walk. 
Indeed, to respect the order of vertices and connected components which is induced by the coupled walks $Z^{\pmb{x},\cdot}$, one has no option but to attach one vertex from the two trees in $\mathcal{F}_1(q-)$ that are merging, and \eqref{Enewedge} means that the parent $L_1(q)$ and the child $L_2(q)$ are picked uniformly from the mass in those components (trees).
\end{remark}

\subsection{Surplus on top of \texorpdfstring{$\mathcal{F}_1$}{F1}}
\label{S:surtopforo}

We refer to any edge of $\mathcal{F}_1(q)$ as a {\em spanning edge arriving before} $q$. 
Note the difference with a similar and weaker definition in terms of $\mathcal{F}_0$. 
Our next goal is to account for the surplus edges in a way which is compatible with the coupled breadth-first walks $Z^{\pmb{x},\cdot}$, and such that the monotonicity of the surplus edges is preserved in time.
\begin{remark}
The construction of the surplus edges in Section \ref{sec:surplus_edges} is not necessarily monotone in time, and it is simpler than the one we are about to describe.
Recall that in the previous setting we only needed to check for the existence of the surplus edges in a certain neighborhood of each vertex, see Figure \ref{F:two} for an illustration of this fact.
However, this is no longer the case here, as our aim is to preserve the monotonicity.
We need to use an external source of randomness evolving in time to check for all the surplus edges among all the available pairs of vertices within a connected component.
\end{remark}

As the reader will see in Section \ref{S:Multi-g}, there is a natural multi-graph that emerges from this construction, and it happens to be the random multi-graph from Bhamidi et al.~\cite{bhamidietal2}, denoted by $(\mathcal{MG}^{\pmb{x}}(t), t\ge 0)$, and already  recalled in the introduction.

\begin{figure}[ht]
\centering
\includegraphics[scale=0.8]{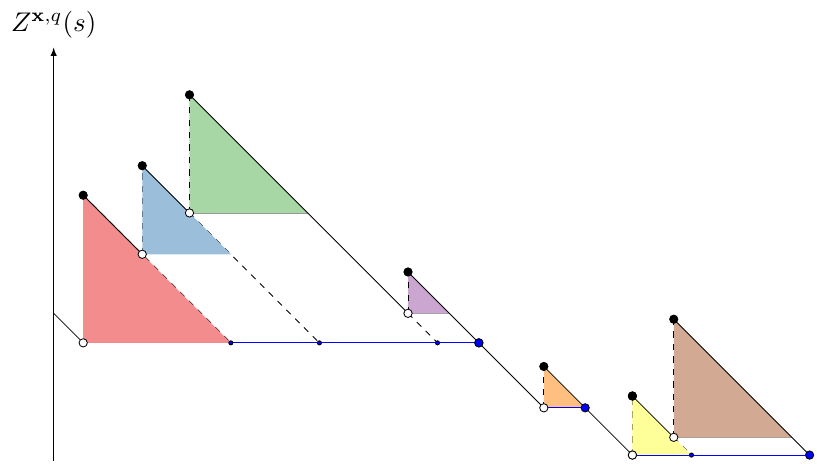}
\caption{The middle excursion is still in its original (and simplest possible) state: it corresponds to the tree consisting of a single (orange) vertex. The excursion to the right corresponds to a tree with two vertices, and the one to the left corresponds to a tree with four vertices. 
The dashed vertical segments are added to indicate the triangles, the $i^{\text{th}}$ of which is matched to vertex $\pi_i$.
The dashed oblique segments indicate zones of interaction which will be described soon.}
\label{F:four}
\end{figure}

Figure \ref{F:four} shows a part of a realization of $Z^{\pmb{x},q}$, with the ``space under the curve'' of the corresponding $B^{\pmb{x},q}$  split into conveniently chosen polygons (to be soon split further into parallelograms), and the triangles marked by different colors for improved clarity.

For $q$ close to $0$ the curve $Z^{\pmb{x},q}$ has only $\operatorname{len}(\pmb{x})$ many (moving) {\em triangular} excursions. As $q$ increases, merging gradually happens, and simultaneously (due to the coupling described in previous sections) the excursions of $Z^{\pmb{x},q}$ become more complex.
It is interesting to describe here
(in Section \ref{sec:ornamented})
the exact structure of this random object, induced by the gradual ``pile up'' of the original triangles (there are $\operatorname{len}(\pmb{x})$ many in total) on ``top of each other''. 

The {\em excursion mosaic} (or simply the {\em mosaic}) is made by drawing at each time $q$ the horizontal (blue in Figure \ref{F:four}) segment on the basis of each excursion of $Z^{\pmb{x},q}$. 
We call these blue segments {\em active baselines}.
Furthermore, if the merging of a pair of trees with roots $\pi_i$ and $\pi_j$ ($i<j$) occurs at time $q$, then
\begin{enumerate}[label=\textbf{O.\arabic*}]
	\item the active baseline starting at $\xi_{(i)}^q$ is extended until the end of the new (just created by this merging event) excursion at time $q$; \label{ornamented_1}
	\item the previously (blue) active baseline starting at $\xi_{(j)}^q$ turns gray, and it is included in the mosaic at any later time $q+z$ as a segment parallel to the active baseline at the vertical level $Z^{\pmb{x},q+z}(\xi_{(j)}^{q+z}-)$; \label{ornamented_2}
	\item for each 
	$l\in [\operatorname{len}(\pmb{x})]$,
	the segment indicating the hypotenuse of the $l^{\text{th}}$ successive triangle is extended in the south-east direction until it meets the active (blue) baseline of its corresponding excursion, so that it connects the points $(\xi_{(l)}^q,Z^{\pmb{x},q}(\xi_{(l)}^q) ) $ and $(\sup(I_l^q), 0))$. We will refer to this extension as the {\em hypotenuse associated to $\pi_l$} (at time $q$). \label{ornamented_3}
\end{enumerate}
We call the gray segments from (\ref{ornamented_2}) the {\em inactive baselines}. If clear from the context, we will simply refer to a blue or gray segment as {\em baseline}.
The $l^{\text{th}}$ successive colored triangle (see Figure \ref{F:four}) is associated to the vertex $\pi_l$.
The parallelogram region 
stretching in the south-east direction from the triangle associated to $\pi_l$ (its top side equals the horizontal cathetus of this triangle) and
situated above the active baseline  (aging see Figure \ref{F:four})
will be called the {\em $\pi_l$-band}.
If $\pi_l$ is the leading vertex of its excursion, the $\pi_l$-band is empty.
We call {\em $\pi_l$-slice}  the union of the $\pi_l$-associated triangle and the $\pi_l$-band.

From now on, we shall call any excursion of $Z^{\pmb{x},q}$, with its corresponding active and inactive baselines (obtained according to the above given  procedure), an  \emph{ornamented excursion.}
The ornamented excursions are intended to record the past history of the component mergers in the evolution of the graph.
Given this information, it will be possible to construct the surplus edges using an external source of randomness, in a way analogous to the construction in Section \ref{sec:surplus_edges}, however adapted to the current setting.

In the next section, we will delve deeper into the ornamented excursions, which is crucial for gaining a comprehensive understanding of these objects. 
However, it is worth noting that the content of the next section may not be directly relevant to the rest of the paper.

\subsubsection{Characterization of the ornamented excursions}
\label{sec:ornamented}

If an ornamented excursion of $Z^{\pmb{x},q}$ traverses $\xi_{(k)}^q$, we say that {\em it carries} $\pi_k$.
It should be clear that the just given mosaic construction and the related definitions can be transposed so that, for each $q$, the reflected process $B^{\pmb{x},q}$ has the same ornamented excursions as $Z^{\pmb{x},q}$, with the difference that in $B^{\pmb{x},q}$ all the excursions start from the abscissa, i.e.\ they are all excursions above level $0$.
 
Note that the ornamented excursions must have some additional ``structural'' properties, due to their construction and the simultaneous breadth-first walk dynamics.
Indeed, almost surely, we have the following conditions
\begin{enumerate}[label=\textbf{O'.\arabic*}]
\item 
\label{ornamented_visual0}
for a fixed $q>0$, no two baselines appear on the same level in the same ornamented excursion,
\item each gray baseline, as constructed by \ref{ornamented_2}, must be a continuous segment (it cannot ``skip'' time intervals), \label{ornamented_visual1}
\item the gray baselines  obey a special kind of monotonicity:  if for $j < l < k$
the baseline at level $Z^{\pmb{x},q}(\xi_{(j)}-)$ reaches the hypotenuse associated to $\pi_l$ (as defined in \ref{ornamented_3}) and
the baseline at level $Z^{\pmb{x},q}(\xi_{(l)}-)$ reaches the hypotenuse associated to $\pi_k$, then 
the baseline at level $Z^{\pmb{x},q}(\xi_{(j)}-)$  must also reach the hypotenuse associated to $\pi_k$.
\label{ornamented_visual2}
\end{enumerate}
Properties \ref{ornamented_visual0}--\ref{ornamented_visual1} are consequences of the way in which the gray baselines are constructed, see \ref{ornamented_2}.
Furthermore, property \ref{ornamented_visual2} is a consequence of the breadth-first walk and the coalescent dynamics: a vertex $\pi_i$ which coalesces at time $q$ with $\pi_j$, such that $j > i$, is immediately connected to all the vertices which were in the connected component of $\pi_j$ just prior to time $q$.

\begin{figure}[h!]
\centering
\begin{subfigure}[b]{0.4\textwidth}
	\centering
	\includegraphics[width=\textwidth]{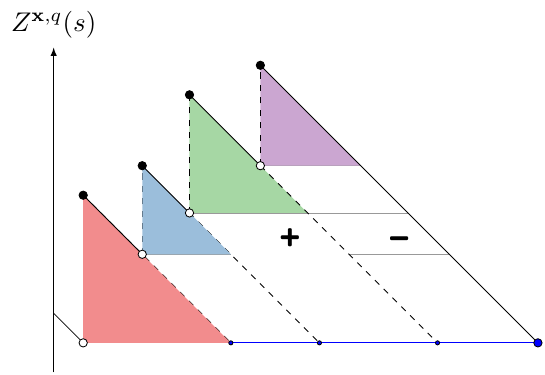}
	\caption{}
	\label{F:six-a}
\end{subfigure}
\hspace{10 pt}
\begin{subfigure}[b]{0.4\textwidth}
	\centering
	\includegraphics[width=\textwidth]{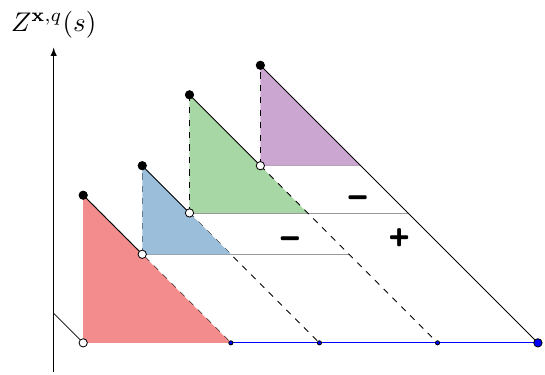}
\caption{}
\label{F:six-b}
\end{subfigure}
\caption{
Two examples of ill-defined ornamented mosaic.
In both cases, in order to obtain a well-defined ornamented excursion, either the gray line(s)
below the minus sign(s) must be removed, or a gray line below the plus sign must be added.
}
\label{F:six}
\end{figure}

Figure \ref{F:six} shows two examples of ill-defined ornamented excursions.
Note that the example in Figure \ref{F:six-a}  (resp.~\ref{F:six-b}) does not satisfy \ref{ornamented_visual1} (resp.~\ref{ornamented_visual2}).

We will next show that the reverse implication is also true: any excursion ornamented according to the three properties just described (the disjoint levels \ref{ornamented_visual0}, the continuity \ref{ornamented_visual1} and the monotonicity \ref{ornamented_visual2}) is a true ornamented excursion, i.e.\ there exists at least one trajectory of the  Aldous' inhomogeneous continuous-time random graph such that the given ornamented excursion appears at some time $q > 0$.

\begin{figure}[h!]
	\begin{subfigure}[b]{0.3\textwidth}
		\centering
		\includegraphics[width=\textwidth]{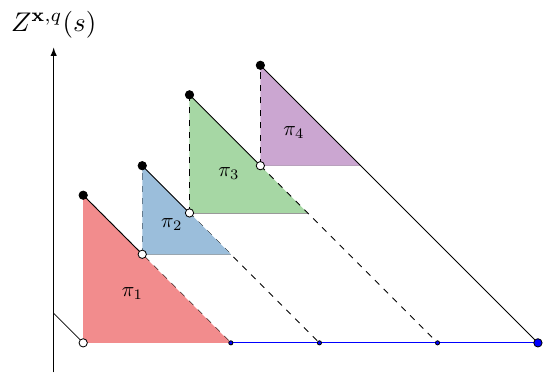}
		\caption{ }
		\label{fig:single_excursion1}
	\end{subfigure}
	\hfill
	\begin{subfigure}[b]{0.3\textwidth}
		\centering
		\includegraphics[width=\textwidth]{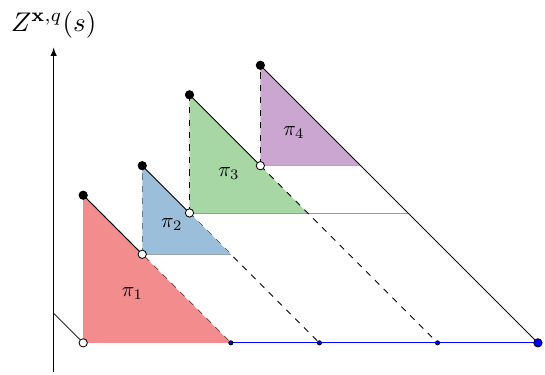}
		\caption{ }
		\label{fig:single_excursion2}
	\end{subfigure}
	\hfill
	\begin{subfigure}[b]{0.3\textwidth}
		\centering
		\includegraphics[width=\textwidth]{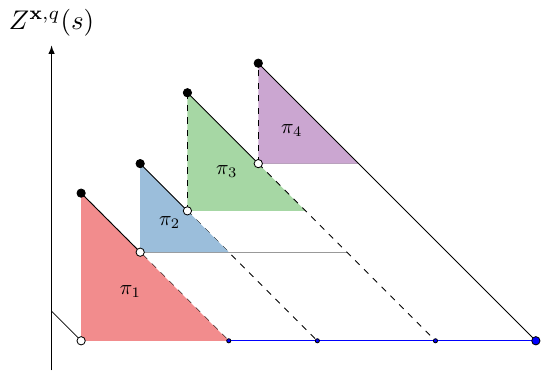}
		\caption{ }
		\label{fig:single_excursion3}
	\end{subfigure}
	\hfill
	\\
	\hfill
	\begin{subfigure}[b]{0.5\textwidth}
		\centering
		\includegraphics[width=0.6\textwidth]{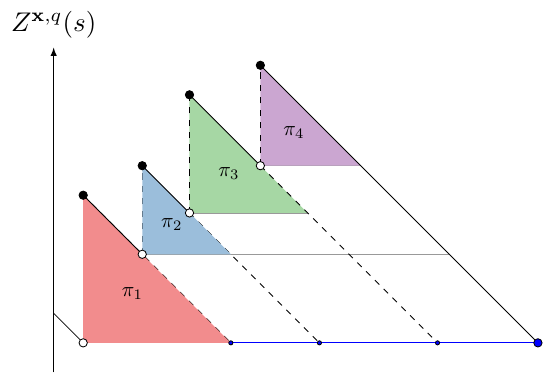}
		\caption{ }
		\label{fig:single_excursion4}
	\end{subfigure}
	\hfill
	\begin{subfigure}[b]{0.5\textwidth}
		\centering
		\includegraphics[width=0.6\textwidth]{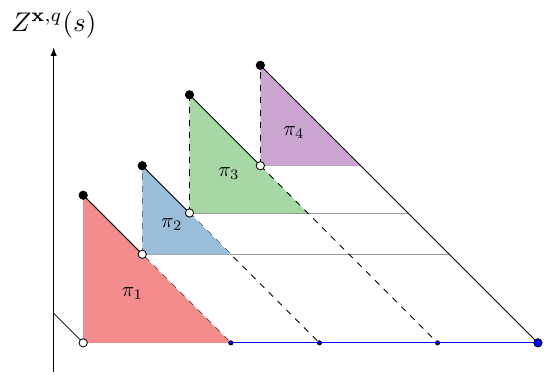}
		\caption{ }
		\label{fig:single_excursion5}
	\end{subfigure}
	\hfill
	\caption{The five possible realizations of the excursion mosaic of four vertices, restricted to a single excursion, viewed at a fixed time $q>0$. 
	}
\label{fig:single_excursion}
\end{figure}

\begin{figure}[h!]
	\begin{subfigure}[b]{0.26\textwidth}
		\centering
		\includegraphics[width=\textwidth]{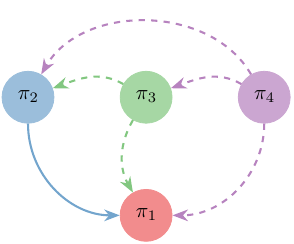}
		\caption{ }
		\label{fig:single_excursion1_tree}
	\end{subfigure}
	\hfill
	\begin{subfigure}[b]{0.16\textwidth}
		\centering
		\includegraphics[width=\textwidth]{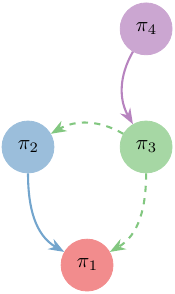}
		\caption{ }
		\label{fig:single_excursion2_tree}
	\end{subfigure}
	\hfill
	\begin{subfigure}[b]{0.16\textwidth}
		\centering
		\includegraphics[width=\textwidth]{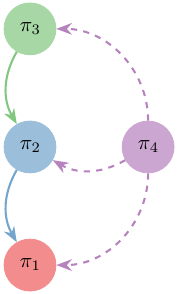}
		\caption{ }
		\label{fig:single_excursion3_tree}
	\end{subfigure}
	\hfill
	\begin{subfigure}[b]{0.16\textwidth}
		\centering
		\includegraphics[width=\textwidth]{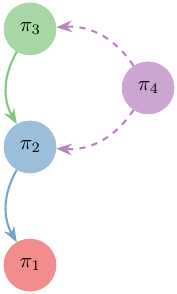}
		\caption{ }
		\label{fig:single_excursion4_tree}
	\end{subfigure}
	\hfill
	\begin{subfigure}[b]{0.055\textwidth}
		\centering
		\includegraphics[width=\textwidth]{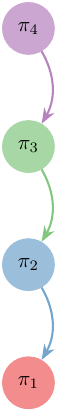}
		\caption{ }
		\label{fig:single_excursion5_tree}
	\end{subfigure}
	\hfill
	\caption{
	Representations of the five trees in $\mathcal{F}_1$ associated to the five ornamented excursions in Figure \ref{fig:single_excursion}.
	A solid arrow indicates an edge in $\mathcal{F}_1$ which is  completely specified by the mosaic.
	The tree corresponding to the fifth image \ref{fig:single_excursion5_tree} (associated to the excursion \ref{fig:single_excursion5}) is completely specified by the mosaic.
	In the other cases, there is at least one edge (indicated by a dashed arrow) in $\mathcal{F}_1$ that needs to be picked via additional randomness. More precisely,
	for each vertex with dashed arrows issued from it, one of the arrows/edges needs to be picked at random according to \eqref{Enewedge}.
	}
	\label{fig:single_excursion_tree}
\end{figure}

If at a given time $q>0$ there is a non-trivial spanning tree of length four in $\mathcal{F}_0$ (or equivalently in $\mathcal{F}_1$), then its corresponding ornamented excursion must have one of the five forms depicted in Figure \ref{fig:single_excursion}.
Figure \ref{fig:single_excursion_tree} shows the possible trees associated to these five ornamented excursions.
The edges completely determined by the ornamented excursions are represented by solid arrows, and those that need to be picked at random according to \eqref{Enewedge} are represented by dashed arrows.

We will illustrate on an example how to construct a trajectory of coalescences in order to obtain a given ornamented excursion.
Without loss of generality, we can suppose that at a given time $q > 0$ the excursion  carries precisely $\{\pi_1, \pi_2, \dots, \pi_k\}$.
The ornamented excursion is then completely determined by the order of coalescence of the vertices $\pi_2, \dots, \pi_k$ with their respective immediately closest leftmost neighbors.
If the leading vertex of an excursion is $\pi_h$, we may write that this excursion is {\em lead by $\pi_h$.}
In the next paragraphs we will describe a construction of a total order on $\{\pi_2, \dots, \pi_k\}$, such that 
when the vertices coalesce according to the enumeration given by this total order,
this results in a given ornamented excursion satisfying \ref{ornamented_visual0}, \ref{ornamented_visual1} and \ref{ornamented_visual2}.

As the first step in our construction, we define
\begin{center}
  a partial order $\succeq$ on the set of vertices carried by an ornamented excursion.  
\end{center}
The leading vertex $\pi_1$ is the greatest element (in symbols, $\pi_1 \succeq \pi_j$, $\forall j$) 
(it is the last element to be visited).
In addition, we say that ${\pi_i} \succeq {\pi_j}$, with $i \le j$, if $i = j$ or if the gray or blue baseline associated to $\pi_i$ (provided it exists, see \ref{ornamented_2}) reaches the hypotenuse line associated to $\pi_j$ (as defined in \ref{ornamented_3}).
This is equivalent to the property that $\pi_j$ was carried by an excursion lead by $\pi_i$ at some earlier time $q'<q$.
Note also that \ref{ornamented_visual2} ensures the transitivity property of the partial order $\succeq$.
Even more is true in view of  \ref{ornamented_visual2}, if ${\pi_i} \succeq {\pi_j}$, with $i < j$, then ${\pi_i} \succeq {\pi_k}$, for every $k$ such that $i \le k \le j$.

\begin{figure}[h]
	\begin{subfigure}[b]{0.24\textwidth}
	\centering
  \includegraphics[width=\textwidth]{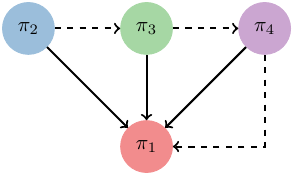}
		\caption{ }
		\label{fig:single_excursion1_diagram}
	\end{subfigure}
	\hfill
	\begin{subfigure}[b]{0.14\textwidth}
		\centering
		\includegraphics[width=\textwidth]{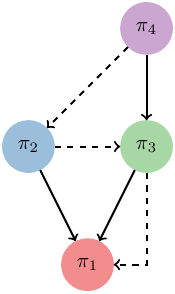}
		\caption{ }
		\label{fig:single_excursion2_diagram}
	\end{subfigure}
	\hfill
	\begin{subfigure}[b]{0.14\textwidth}
		\centering
		\includegraphics[width=\textwidth]{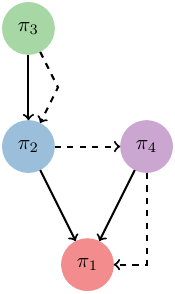}
		\caption{ }
		\label{fig:single_excursion3_diagram}
	\end{subfigure}
	\hfill
	\begin{subfigure}[b]{0.14\textwidth}
		\centering
		\includegraphics[width=\textwidth]{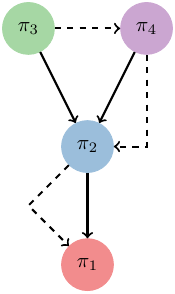}
		\caption{ }
		\label{fig:single_excursion4_diagram}
	\end{subfigure}
	\hfill
	\begin{subfigure}[b]{0.079\textwidth}
		\centering
		\includegraphics[width=\textwidth]{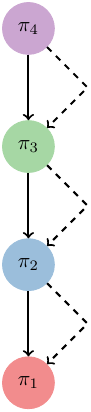}
		\caption{ }
		\label{fig:single_excursion5_diagram}
	\end{subfigure}
	\hfill
	\caption{
	The tree formed by the solid (resp.~dashed) arrows is the Hasse diagram associated to the partial order $\succeq$ (resp.~total order $\succeq^\star$).
	Collapsing the nodes according to the paths defined by the dashed arrows will produce the ornamented excursions in Figure \ref{fig:single_excursion}.	
	}
	\label{fig:single_excursion_diagram}
\end{figure}

\smallskip 
Let us recall that the \textit{Hasse diagram} associated to the partial order $\succeq$ is the directed graph with edges $\pi_i \leftarrow  \pi_j$ whenever $\pi_i \succeq \pi_j$ and there is no other $\pi_k$ such that $\pi_i \succeq \pi_k \succeq \pi_j$ (in this case it is said that $\pi_i$ \textit{covers} $\pi_j$ in the logic literature language).
The {Hasse diagram} associated to $\succeq$ is consequently a rooted tree whose root is the first element in the excursion, i.e.\ the greatest element according to $\succeq$.
The trees given by the solid arrows in Figure \ref{fig:single_excursion_diagram} are the Hasse diagrams associated to the partial order $\succeq$ obtained from the excursions in Figure \ref{fig:single_excursion}.

Furthermore, we will use
\begin{center}
a total order $\succeq^\star$ on the set of vertices carried by an ornamented excursion. \end{center}
 which we define next:
if $\pi_i \succeq \pi_j$ then $\pi_i \succeq^\star \pi_j$, and  
 if $\pi_i$ and $\pi_j$ are not ordered according to $\succeq$ but they both belong to the same 
 generation of the Hasse diagram, then $i\geq j$ implies $\pi_i \succeq^\star \pi_j$.
It is easy to check that $\succeq^\star$ is indeed a total order.

We refer the reader to Figure \ref{fig:single_excursion_diagram} and Figure \ref{F:seven}.
Figure \ref{F:seven-a} shows a more complicated ornamented excursion carrying eight vertices.
Figure \ref{F:seven-b} shows the associated tree where, similarly to Figure \ref{fig:single_excursion_tree}, the edges completely determined by the excursion are represented by solid arrows and those that need to be picked according to \eqref{Enewedge} are represented by dashed arrows.
In addition, Figure \ref{F:seven-c} shows the respective Hasse diagrams associated to the orders $\succeq$ and $\succeq^\star$ defined on the set of vertices (excluding the root), with solid and dashed arrows, respectively.

Each merger event concerns two successive excursions, such that their baselines touch at the coalescent time. The excursion to the left (resp.~excursion to the right) is lead by some $\pi_{l}$ (resp.~$\pi_{r}$) where it must be $l<r$.
An arrow in the Hasse diagram (solid lines in Figure  \ref{fig:single_excursion_diagram}) between $\pi_r $ and $\pi_l $ where $l<r$ indicates that the excursion lead by $\pi_r $ is to be merged with an excursion lead by $\pi_l$ (and necessarily carrying $\{\pi_l,\ldots, \pi_{r-1}\})$. 
With the convention that the edges in the Hasse diagram point from children to parents, let us
denote by ${\rm par}(v)$ the parent of $v$. 
The total order $\succeq^\star$ (dashed lines in Figure \ref{fig:single_excursion_diagram} and Figure \ref{F:seven-c}) specifies the succession of merger events.
More precisely, let $\tau$ be the permutation of $2,\ldots,k$ such that $\pi_1 \succeq^\star \pi_{\tau_2} \succeq^\star \pi_{\tau_3}\succeq^\star \ldots \succeq^\star \pi_{\tau_k}$.
At step $m=0$ our configuration consists of $k$ isolated excursions (or vertices) and for each $m\in \{1,\ldots,k-1\}$ at step $m$ of the construction we take the configuration from step $m-1$ and merge the excursion lead by $\pi_{\tau_{k-m+1}}$ with the excursion lead by ${\rm par}(\pi_{\tau_{k-m+1}})$, resulting in the configuration at step $m$.
Due to the above reasoning, the coalescence of $\pi_{\tau_{k-m+1}}$ with ${\rm par}(\pi_{\tau_{k-m+1}})$ produces exactly the gray baseline associated to 
$\pi_{\tau_{k-m+1}}$
in the given ornamented excursion.

Note that coalescing the vertices in a different order than the one induced by $\succeq^\star$ could produce undesired gray baselines. 
Indeed, for two vertices such that $\pi_i \succeq^\star \pi_k$ and $i > k$, the coalescence of vertex $\pi_i$ with ${\rm par}(\pi_i)$  before the coalescence of $\pi_k$ with ${\rm par}(\pi_k)$  could produce a gray baseline associated to $\pi_k$ which is longer than the one given by the ornamented excursion.
See, for instance, the difference between the Hasse diagrams in Figures \ref{fig:single_excursion4_diagram} and \ref{fig:single_excursion5_diagram}, related to the ornamented excursions in Figures \ref{fig:single_excursion4} and \ref{fig:single_excursion5}, respectively, focusing on $i = 4$ and $k = 3$.

The total order induced by $\succeq^\star$ (dashed arrows) associated to the ornamented excursions in Figure \ref{fig:single_excursion} and represented by the Hasse diagrams in Figure \ref{fig:single_excursion_diagram} are summarized in Table \ref{table:total_ordering}.
Notice that, even though the total order defined for a given ornamented excursion is a natural choice, it is not the only ordering of coalescence of vertices resulting in the given ornamented excursion.
It can be easily checked that the total order $\pi_1 \succeq^\star \pi_3 \succeq^\star \pi_4 \succeq^\star \pi_2$, which is not listed in Table \ref{table:total_ordering}, also produces the ornamented excursion in Figure \ref{fig:single_excursion2}, i.e.\ the same produced by the ordering $\pi_1 \succeq^\star \pi_3 \succeq^\star \pi_2 \succeq^\star \pi_4$.

\begin{table}[h!]
	\centering
\begin{tabular}{|c|c|c|}
	\hline
	\rule[-1ex]{0pt}{2.5ex} Ornamented excursion & Hasse diagram & Total order \\
	\hline
	\rule[-1ex]{0pt}{2.5ex} Figure  \ref{fig:single_excursion1} & Figure  \ref{fig:single_excursion1_diagram} & $\pi_1 \succeq^\star \pi_4 \succeq^\star \pi_3 \succeq^\star \pi_2$ \\
	\hline
	\rule[-1ex]{0pt}{2.5ex} Figure  \ref{fig:single_excursion2} & Figure  \ref{fig:single_excursion2_diagram} & $\pi_1 \succeq^\star \pi_3 \succeq^\star \pi_2 \succeq^\star \pi_4$ \\
	\hline
	\rule[-1ex]{0pt}{2.5ex} Figure  \ref{fig:single_excursion3} & Figure  \ref{fig:single_excursion3_diagram} & $\pi_1 \succeq^\star \pi_4 \succeq^\star \pi_2 \succeq^\star \pi_3$ \\
	\hline
	\rule[-1ex]{0pt}{2.5ex} Figure  \ref{fig:single_excursion4} & Figure  \ref{fig:single_excursion4_diagram} & $\pi_1 \succeq^\star \pi_2 \succeq^\star \pi_4 \succeq^\star \pi_3$ \\
	\hline
	\rule[-1ex]{0pt}{2.5ex} Figure  \ref{fig:single_excursion5} & Figure  \ref{fig:single_excursion5_diagram} & $\pi_1 \succeq^\star \pi_2 \succeq^\star \pi_3 \succeq^\star \pi_4$ \\
	\hline
\end{tabular}
\caption{
	Total order induced by the ornamented excursions in Figure \ref{fig:single_excursion}, whose Hasse diagrams are represented by dashed arrows in Figure \ref{fig:single_excursion_diagram}.
}
\label{table:total_ordering}
\end{table}
%

\vspace{-0.5cm}

\subsubsection{The surplus edges dynamic}
\label{S:surplus_edge_dynamic}
\begin{figure}[h]
	\centering
	\includegraphics[width=0.5\linewidth]{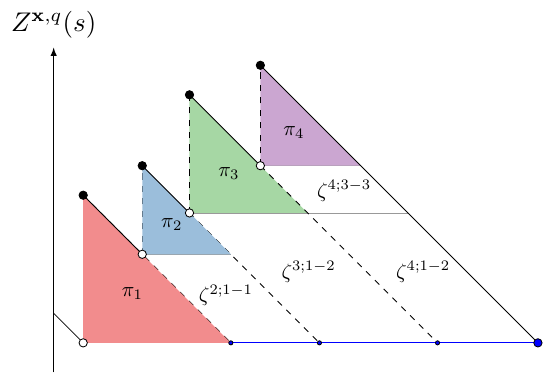}
	\includegraphics[width=0.2\linewidth]{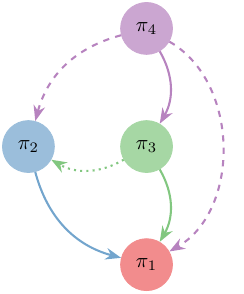}
	\caption{The ornamented excursion is the same as in Figure \ref{fig:single_excursion2}. 
	In the tree, it is assumed that the edges $\pi_2 \rightarrow \pi_1$, $\pi_3 \rightarrow \pi_1$,  and  $\pi_4 \rightarrow \pi_3$ were sampled according to \eqref{Enewedge} and thus included in $\mathcal{F}_1$. 
	The four independent PPP $\zeta^{2;1-1}$, $\zeta^{3;1-2}$, $\zeta^{4;1-2}$ and $\zeta^{4;3-3}$, with rates $x_{\pi_2} \cdot x_{\pi_1}$, $x_{\pi_3} (x_{\pi_2} + x_{\pi_1})$, $x_{\pi_4} (x_{\pi_2} + x_{\pi_1})$ and $x_{\pi_4} \cdot x_{\pi_3}$, respectively, are matched to the four parallelogram regions as shown in the excursion.
	The dotted green arrow represents the surplus edge that will be created using $\zeta^{3;1-2}$ and the dashed violet arrows represent the surplus edges that will be created using $\zeta^{4;1-2}$.
	}
	\label{fig:singleexcursiontree2-xi}
\end{figure}

Let us first show  on an example how the surplus edges can be superimposed in a consistent way.
Consider again the example from Figure \ref{fig:single_excursion2}. 
Let us enumerate the four parallelograms specified by the mosaic in some way, for example traverse them row-by-row from left to right, as an analogue to the breadth-first order (see Figure \ref{F:one}). 
Let $\zeta^{2;1-1}$, $\zeta^{3;1-2}$, $\zeta^{4;1-2}$ and $\zeta^{4;3-3}$ be independent Poisson point processes (PPP), matched respectively to these four regions, as shown in Figure \ref{fig:singleexcursiontree2-xi}.
Note that the the notation $\zeta^{l; j-k}$ with $j < k < l$ means that the parallelogram is contained in the $\pi_l$-band and its southern (resp.\ northern) base is contained in the baseline associated to $\pi_j$ (resp.\ $\pi_{k + 1}$).

As $q$ increases, the base of each parallelogram stays fixed in length (although it moves closer to the origin), while its height increases. 
Let a point arrive to $\zeta^{2;1-1}$ at rate $x_{\pi_2} \cdot x_{\pi_1}$, to $\zeta^{3;1-2}$ at rate $x_{\pi_3} (x_{\pi_2} + x_{\pi_1})$, to $\zeta^{4;1-2}$ at rate $x_{\pi_4} (x_{\pi_2} + x_{\pi_1})$, and to $\zeta^{4;3-3}$ at rate $x_{\pi_4} \cdot x_{\pi_3}$.
As it turns out, most of these Poisson point processes will not be needed now, but they are nevertheless defined with an intention of later use.
Let us assume that at a given time $q$, due to randomization \eqref{Enewedge}, the realization of the corresponding random tree has the following three edges:  $\pi_2 \rightarrow \pi_1$, $\pi_3 \rightarrow \pi_1$,  and  $\pi_4 \rightarrow \pi_3$, indicated by solid lines in Figure \ref{fig:singleexcursiontree2-xi}.
We only need to watch for marks in $\zeta^{3;1-2}$ and $\zeta^{4;1-2}$.
When a new mark arrives to $\zeta^{3;1-2}$, 
with probability $x_{\pi_2}/(x_{\pi_1}+x_{\pi_2})$ the process creates an edge $\pi_3 \rightarrow \pi_2$, unless this edge already exists.
This surplus edge is represented by a dotted arrow in Figure \ref{fig:singleexcursiontree2-xi}. 
Nothing happens with the remaining probability.
When a new mark arrives to $\zeta^{4;1-2}$,
with probability $x_{\pi_2}/(x_{\pi_1}+x_{\pi_2})$ (resp.~$x_{\pi_1}/(x_{\pi_1}+x_{\pi_2})$) the process creates an edge $\pi_4 \rightarrow \pi_2$ (resp.~$\pi_4 \rightarrow \pi_1$),  unless it already exists.
These surplus edges are represented by dashed arrows in Figure \ref{fig:singleexcursiontree2-xi}. 
It is likely clear to the reader that these transitions are chosen with the purpose of preserving the random graph transitions within the connected components.

\smallskip
In the general case, 
 one takes the collection of independent Poisson processes $\zeta^{l;j-k}$, where $j,k,l$ run over all the indices in $[\operatorname{len}(\pmb{x})]$ such that
\begin{equation*}
\label{Econstr} 
j\leq k\leq l, \mbox{ and if }k=l\mbox{ then also }j=l.
\end{equation*} 
The processes $(\zeta^{l;l-l})_l$ play a special role, to be explained in Section \ref{S:Multi-g}. 
If $l>k$, then $\zeta^{l;j-k}$ is in charge of generating a surplus edge from $\pi_l$ to some vertex in the range $\{\pi_j,\ldots,\pi_k\}$, but {\bf only when compatible} with the excursion mosaic. 
More precisely, $\zeta^{l;j-k}$ will be {\em active} starting from  the time $T_{l;j-k}$ at which the ornamented excursion of $Z^{\pmb{x},\cdot}$ carrying $\pi_l$ (as its root or a non-root vertex) merges with an ornamented excursion carrying precisely $\{\pi_j,\ldots,\pi_k\}$. 
In this way, the random time $T_{l;j-k}$ depends on the mosaic, and it is finite if and only if the just described merging configuration occurs. 
It can happen that $j=k<l$, in the case where the ornamented excursion carrying $\pi_l$ (with leading vertex $\pi_{j+1}$) merges with a simple triangular excursion carrying only $\pi_j$.
On the event $\{T_{l;j-k} = \infty\}$, the corresponding $\zeta^{l;j-k}$ is never activated.
On $\{T_{l;j-k} < \infty\}$, the behavior of $\zeta^{l;j-k}$ after $T_{l;j-k}$ is a generalization of the one given in the example above (see e.g.~the definition of $\zeta^{3;1-2}$).
More precisely, the points arrive to $\zeta^{l;j-k}$ at rate $x_{\pi_l} (x_{\pi_j} + \ldots + x_{\pi_k})$. 
The first (spanning) edge is created according to \eqref{Enewedge} when the merging occurs.
The surplus edges are created as follows: given an arrival to $\zeta^{l;j-k}$ at time $q > T_{l;j-k}$, an index $I$ is drawn (independently of the past of the mosaic, of $\mathcal{F}_1$, and of the surplus edge data), so that $I=i$ with probability $x_{\pi_i}/( x_{\pi_j} + \ldots + x_{\pi_k})$ for each $i \in \{j,\ldots,k\}$. Given $I$, the surplus edge $l\rightarrow I$ is created, unless it already exists. 
For each $q>0$, call any edge created in this way before time $q$ a {\em surplus edge arriving before} $q$.
Let $G_1^{\pmb{x}}(q)=(V,E_1(q))$, where $E_1(q)$ is the union of the spanning edges in $\mathcal{F}_0(q)$ and surplus edges arriving before $q$.
Then it is clear that $G_1^{\pmb{x}}$ is a monotone random graph process: $G_1^{\pmb{x}}(q_1) \subset G_1^{\pmb{x}}(q_2)$ whenever $q_1\leq q_2$. In fact, we can derive a stronger claim, which may serve as the main motivation for the additional multi-graph construction presented in the next section.

The following is an immediate and important conclusion of the observations given above.
\begin{theorem}[Encoding of the simple graph process]
\label{thm:strictatqfull}
For any general initial positive weights $\pmb{x}$,
the just constructed graph-valued process
$(G_1^{\pmb{x}}(q),\,q\geq 0)$
is a realization of the Aldous' (inhomogeneous) continuous-time random graph $(\mathcal{G}^{\pmb{x}} (q), q \ge 0)$.
\end{theorem}
 
\subsection{Proof of Theorem \ref{thm:encodingAMC}}
\label{S:Multi-g}
Here we focus on the second construction above (see Section \ref{S:surtopforo}).
In particular, recall the excursion mosaic, and the family of compatible Poisson point processes
$\zeta^{l;j-k}$, where $j,k,l$ satisfy the constraints given in \eqref{Econstr}.
For each $l$, $\zeta^{l;l-l}$ should now be matched at time $q$ to the triangle spanned by the points 
$(\xi_{(l)}^q, Z^{\pmb{x},q}(\xi_{(l)}^q-))$, $(\xi_{(l)}^q, Z^{\pmb{x},q}(\xi_{(l)}^q))$
and
$(\xi_{(l)}^q+x_{\pi_l}, Z^{\pmb{x},q}(\xi_{(l)}^q-))$ (or equivalently, to that spanned by $(\xi_{(l)}^q, B^{\pmb{x},q}(\xi_{(l)}^q-))$,  $(\xi_{(l)}^q, B^{\pmb{x},q}(\xi_{(l)}^q))$ and $(\xi_{(l)}^q+ x_{\pi_l}, B^{\pmb{x},q}(\xi_{(l)}^q-))$).
These are the colored triangles in the mosaic illustrated in Figure \ref{fig:singleexcursiontree2-xi}.
The process $\zeta^{l;l-l}$ is active already at time $0$ (the triangular excursions exist from the very beginning). 
We therefore define $T_{l;l-l}\equiv 0$ almost surely, for each $l\in [\operatorname{len}(\pmb{x})]$. 
Points arrive to $\zeta^{l;l-l}$ at rate $x_{\pi_l}^2/2$, and at the time of each arrival, a self-loop $\pi_l \rightarrow \pi_l$ is created.
\begin{remark}
\label{R:isolated_vertex}
The factor of $1/2$ is natural if one thinks of each original block as continuous ``spread'' of mass, and of each self-loop as an edge between two different points in this block. 
If $x$ is discretized into $n$ subblocks of equal mass $x/n$ for some large $n$, and if an edge between the $i^{\text{th}}$ and the $j^{\text{th}}$ subblock arrives at the multiplicative rate $(x/n)^2$, then the total rate of edge arrivals equals $(x/n)^2 \cdot {n \choose 2} = \frac{x^2}{2} (1- 1/n) \approx x^2/2$.
Not surprisingly, this rate is also the area of the triangle to which $\zeta^{l; l-l}$ is matched.
\end{remark}

As in Section \ref{S:surplus_edge_dynamic}, for $l>k$ the counting process $\zeta^{l;j-k}$ is activated at time $T_{l;j-k}$, hence never on  $\{T_{l;j-k}=\infty\}$. 
After activation, the points arrive to $\zeta^{l;j-k}$ again at rate  $x_{\pi_l} (x_{\pi_j} + \ldots + x_{\pi_k})$.
The surplus multi-edges are created as before, but without a ``uniqueness'' constraint: given an arrival to $\zeta^{l;j-k}$ at time $q$, an index $I$ is drawn in the same way as before, and a new surplus (multi-)edge $l\rightarrow I$ is created at time $q$. 
Let us call this process $(MG_1^{\pmb{x}}(q),\,q\geq 0)$.
\begin{theorem}[Encoding of the multi-graph process]
\label{Theorem3.2B}
For any general initial positive weights $\pmb{x}$,
the just constructed multi-graph-valued process
$(MG_1^{\pmb{x}}(q),\,q\geq 0)$
is a realization of the continuous-time multi-graph-valued process $(\mathcal{MG}^{\pmb{x}} (q), q \ge 0)$.
\end{theorem}
We already explained how $\zeta^{l;l-l}$ can be matched to the $l^{\text{th}}$ triangular region under the curve $B^{\pmb{x},\cdot}$, or equivalently, the triangle associated with $\pi_l$.
It is useful to point out here that $\zeta^{l;j-k}$ (on $\{T_{l;j-k} < \infty\}$) can analogously be matched to a parallelogram shaped region (evolving in time) on the ornamented mosaic, for any choice  of $l>k\geq j$ such that $\zeta^{l;j-k}$ gets activated. 
Before time $T_{l;j-k}$ this parallelogram does not exist, exactly at time $T_{l;j-k}$ it has height $0$, and its height (strictly) increases at any future time. 
Indeed, this parallelogram of constant base length $x_{\pi_l}$ is created at time $T_{l;j-k}$ by the excursion mosaic, and at any time $z>T_{l;j-k}$ it is the bounded region specified by the four lines
\begin{equation}
\label{Elines}
\begin{array}{ll}
y= - \left(x - \xi_{(l)}^z\right) +Z^{\pmb{x},z}\left(\xi_{(l)}^{z}-\right) & \ \mbox{\small{``left'' boundary}},\\
y= Z^{\pmb{x},z}\left(\xi_{(j)}^{z}-\right) & \ \mbox{\small{``bottom'' boundary}},\\
y= - \left(x - \xi_{(l)}^z\right) +Z^{\pmb{x},z}\left(\xi_{(l)}^{z}\right) & \ \mbox{\small{``right'' boundary}},\\
y= Z^{\pmb{x},z}\left(\xi_{(k+1)}^{z}-\right) & \ \mbox{\small{``top'' boundary}}.\\
\end{array}
\end{equation}

\begin{remark}\label{remark:parallelogram}
It is easily derived from the construction that 
on $\{T_{l;j-k}<\infty\}$,  the excursion carrying $\pi_l$ must have $\pi_{k+1}$ as its leading vertex (or equivalently, $\pi_{k+1}$ is the root of the corresponding spanning tree) just prior to its merger with the excursion carrying $\{\pi_j, \ldots,\pi_k \} $.
\end{remark}

\begin{figure}[h]
\centering
\begin{subfigure}[b]{0.5\textwidth}
	\centering
	\includegraphics[width=\textwidth]{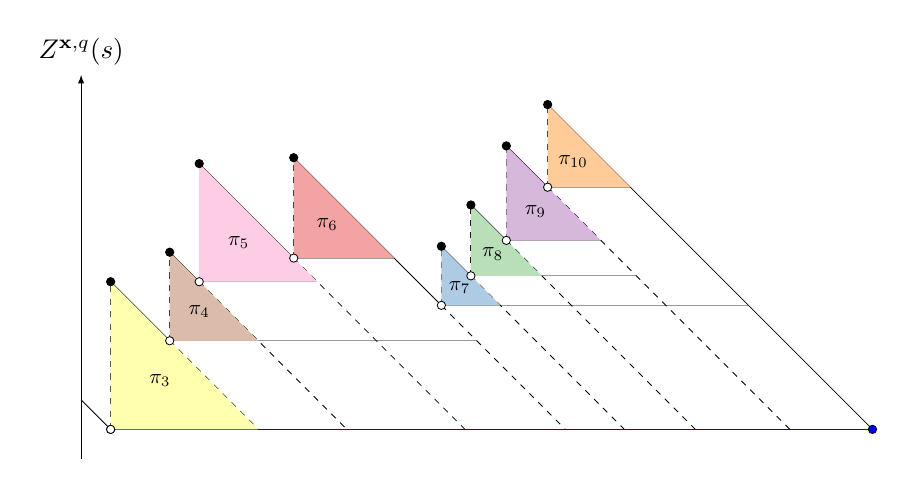}
	\caption{ }
	\label{F:seven-a}
\end{subfigure}
\hfill
\begin{subfigure}[b]{0.2\textwidth}
	\centering
	\includegraphics[width=\textwidth]{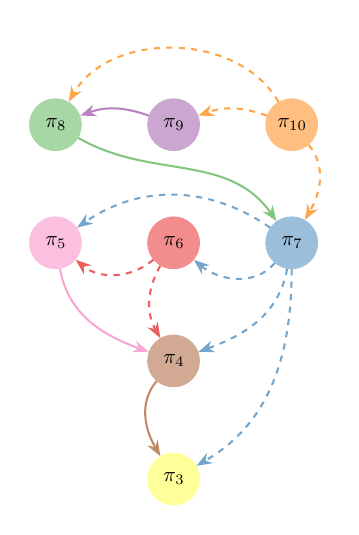}
	\caption{ }
	\label{F:seven-b}
\end{subfigure}
\hfill
\begin{subfigure}[b]{0.265\textwidth}
	\centering
	\includegraphics[width=\textwidth]{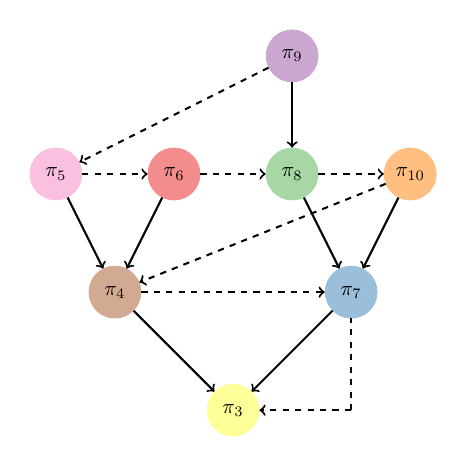}
	\caption{ }
	\label{F:seven-c}
\end{subfigure}
\hfill
\caption{
	The tree \subref{F:seven-b} associated to the ornamented excursion \subref{F:seven-a} is partially determined by the mosaic.
	The dashed edges in \subref{F:seven-b} indicate various options, to be decided according to \eqref{Enewedge}.
	The solid arrows in \subref{F:seven-c} are the edges in the Hasse diagram related to the order $\succeq$ associated to the ornamented excursion \subref{F:seven-a}, as described in Section \ref{sec:ornamented}.
	The dashed arrows in \subref{F:seven-c} represent the Hasse diagram of the total order $\succeq^\star$, also defined in Section \ref{sec:ornamented}.
	For the violet vertex $\pi_9$, there are three active Poisson point processes $\zeta^{9;8-8}$, $\zeta^{9;7-7}$ and $\zeta^{9;3-6}$ running at that time.
	For the orange vertex $\pi_{10}$, there are only two active Poisson point processes $\zeta^{10;7-9}$, and  $\zeta^{10;3-6}$ running.
	The process $\zeta^{9;7-7}$ corresponds to the middle parallelogram in the $\pi_9$-band, containing the triangle marked in violet, and $\zeta^{10;3-6}$ corresponds to the lower parallelogram in the $\pi_{10}$-band, containing the triangle marked in orange.
There is no $\zeta^{7;3-5}$, why? Find $\zeta^{6;3-3}$ and $\zeta^{6;4-5}$.
}
\label{F:seven}
\end{figure}

Figure \ref{F:seven} shows how the mosaic drawn in previous figures might look at a later time. 
For the sake of illustration, let us assume that the eight jumps (vertices) in Figure \ref{F:seven} 
are $\pi_3,\ldots,\pi_{10}$, where $\pi_3$ is indicated in yellow, and $\pi_{10}$ in orange (the jumps at $\xi_{(1)}^q$ and at $\xi_{(2)}^q$ are not contained in the plot region).

The excursion mosaic is an interesting object due to the following claim in particular.
\begin{proposition}
\label{P:rates}
For any choice of $j,k,l$
satisfying \eqref{Econstr},
and for any fixed $q$, on $\{T_{l;j-k}\leq q\}$,
the cumulative arrival rate to $\zeta^{l;j-k}$ up to time $q$ equals almost surely the area of the region in the excursion mosaic matched to $\zeta^{l;j-k}$ at time $q$, multiplied by $q$.
\end{proposition}
\begin{proof}
If $j=k=l$, the statement is trivial to check (the region is the right isosceles triangle, see Remark \ref{R:isolated_vertex}).
Suppose that $l>k$, so that the corresponding region in the mosaic is a parallelogram. 
Then clearly the total arrival rate to $\zeta^{l;j-k}$ before time $q$ equals $(q-T_{l;j-k})^+ \cdot x_{\pi_l} (x_{\pi_j} + \ldots + x_{\pi_k})$.

Abbreviate 
$$\chi := x_{\pi_j} + \ldots + x_{\pi_k}.$$
Then on $\{T_{l;j-k} \geq q \}$ the total rate above can be rewritten as 
$$
q \left(\chi - \chi \cdot \frac{ T_{l;j-k}}{q}\right) \cdot x_{\pi_l}.
$$ 
The claim we need to prove is therefore 
equivalent to the fact that the height of the parallelogram specified in \eqref{Elines} equals the term in parentheses (see also Remark \ref{remark:parallelogram}). 
At time $q$, the height of the parallelogram from \eqref{Elines} is 
\begin{equation}\label{eq:height_parallelogram1}
    Z^{\pmb{x},q}\left( \xi_{(k+1)}^{q}- \right) - Z^{\pmb{x},q} \left(\xi_{(j)}^{q}- \right) \equiv \chi - \frac{\xi_{(k+1)} - \xi_{(j)}}{q}.
\end{equation}
It was already noted that $T_{l;j-k}$ (whenever finite) is the time of merging of two excursions (one carrying $\pi_l$ and the other carrying exactly $\pi_j,\ldots,\pi_k$ immediately prior to $T_{l;j-k}$). 
Recalling Remark \ref{remark:parallelogram} again, and using the definition of $(Z^{\pmb{x},q})_{q>0}$, we conclude that the time $T_{l;j-k}$ of  merging of these two excursions must satisfy the identity
\begin{equation}\label{eq:height_parallelogram2}
    \xi_{(k+1)} - \xi_{(j)} = \chi \cdot T_{l;j-k}.
\end{equation}
Plugging \eqref{eq:height_parallelogram2} into \eqref{eq:height_parallelogram1}, it is  simple to verify that on $\{T_{l;j-k}\leq q\}$ the quantity $Z^{\pmb{x},q} \left( \xi_{(k+1)}^{q}- \right) - Z^{\pmb{x},q} \left( \xi_{(j)}^{q}- \right)$
equals $\chi \cdot (1 - T_{l;j-k}/{q})$, almost surely, 
which is the required identity.
\end{proof}

Recall the 
slices defined from the graph of the the non-reflected process $Z^{\pmb{x},z}$ in Section \ref{S:surtopforo}.
We could equivalently consider the paths of $B^{\pmb{x},z}$.
Notice that the $\pi_l$-slice (under the curve of $B^{\pmb{x},z}$) is the bounded (random) region specified by the lines 
\begin{equation}
\begin{array}{ll}
y= - \left( x - \xi_{(l)}^q \right) + B^{\pmb{x},q} \left( \xi_{(l)}^{q}- \right) & \ \mbox{\small{``left lower'' boundary}},\\
x= \xi_{(l)}^{q}, & \ \mbox{\small{``left upper'' boundary}}\\
y= - \left( x - \xi_{(l)}^q \right) + B^{\pmb{x},q} \left( \xi_{(l)}^{q} \right) & \ \mbox{\small{``right'' boundary}},\\
y= 0 & \ \mbox{\small{``bottom'' boundary}},\\
\end{array}
\nonumber
\end{equation}
We recall that the $\pi_l$-slice is split by the mosaic into the right isosceles triangle of area $(x_{\pi_l})^2/2$ and (possibly) additional parallelograms, each of which corresponds to $\zeta^{l;j-k}$ for some $j\leq k<l$ such that $T_{l;j-k}\leq q$.
In this way, the region under the excursion of $B^{\pmb{x},z}$ carrying exactly $\pi_h,\pi_{h+1},\ldots, \pi_{h+r}$ at time $q$ is the union of the $\pi_h$-slice, the $\pi_{h+1}$-slice, $\ldots$, and the $\pi_{h+r}$-slice.
The intersections of adjacent trinales, parallelograms, bands and slices are sets (segments) of zero Lebesgue measure.
Thus, the cumulative rate of (oriented) surplus edges issued from $\pi_l$ before time $q$ is the area of $\pi_l$-slice at time $q$, multiplied by $q$.
Summing over all the slices in the excursion carrying $\pi_h,\pi_{h+1},\ldots, \pi_{h+r}$ at time $q$ leads to the next result.
\begin{lemma}[Area below the excursion and surplus edges] \label{coro:area_below_B}
    At time $q$ the conditional cumulative rate of surplus edges in the component of $\mathcal{MG}^{\pmb{x}}(q)$ with vertices $\pi_h,\pi_{h+1},$ $\ldots, \pi_{h+r}$ is the area of the excursion of $B^{\pmb{x},q}$ carrying exactly $\pi_h,\pi_{h+1},\ldots, \pi_{h+r}$, multiplied by $q$. Moreover, the surplus edge counts are conditionally independent over different connected components.
\end{lemma}

\begin{proof}[Proof of Theorem \ref{thm:encodingAMC}]
Let $\mathcal{C}_i$ be the $i^{\text{th}}$ largest component of $\mathcal{MG}^{\pmb{x}}(q)$, as constructed by Theorem \ref{Theorem3.2B}\,(B).
Lemma \ref{coro:area_below_B}, and elementary facts on the sum of independent Poisson random variables, implies that the AMC corresponding to $\mathcal{MG}^{\pmb{x}}(q)$ has the same distribution as $(\pmb{X}^{\pmb{x}}, \pmb{N}^{\pmb{x}})$ coupled to $B^{\pmb{x},q}$ via \eqref{eq:def_encoding}.
\end{proof}

\section{The sBFW  encodes the AMC dynamic}
\label{sec:sBFWencodesAMCdynamic}
In this section we prove Theorem \ref{thm:dynamicAMC2}.
Recall from the Introduction that $\rho_s$, the total length of load-free periods of $B^{\pmb{x},q}$ up to time $s$, satisfies
\[
    \rho_s = \int_0^s \mathbb{1}_{\{B^{\pmb{x},q}(u) = 0\}} \mathrm{d}u = - \inf_{u \le s} Z^{\pmb{x},q}(u).
\]
Also recall 
the process $\mathfrak{B}^{\pmb{x}, q}$ defined in
\eqref{eq:def_B_curly}, as well as its 
alternative expression on $\left[ 0, \sum_{i = 1}^n x_i \right]$,
\begin{equation}
\label{E:formula_curly_B}
\mathfrak{B}^{\pmb{x},q}(s) = \sum_{i = 1}^n x_{i} \cdot \mathbb{1}_{\{ \Xi(\xi_i/q) \le s \}} - s, \;\; \text{ where }\ \Xi: s \mapsto s - \rho_{s} = \int_0^s \mathbb{1}_{\{ B^{\pmb{x}, q}(u) > 0 \}} \mathrm{d} u.  
\end{equation}
The finite sequences of excursions (above $0$) of $\mathfrak{B}^{\pmb{x},q}$ and of $B^{\pmb{x},q}$ are the same.
In particular, the lengths of excursions, when listed chronologically with respect to either $\mathfrak{B}^{\pmb{x},q}$ or $B^{\pmb{x},q}$, are the same.
Due to the definitions, the $i^{\text{th}}$ largest excursion of $\mathfrak{B}^{\pmb{x},q}$ occurs  during the interval $\left[ \Xi(a_i(q)) , \Xi(a_i(q)) + X_i(q)  \right]$, where $a_i(q)$ is 
the starting time of the $i^{\text{th}}$ longest excursion of $B^{\pmb{x},q}$ (its length is $X_i(q)$),  
and where $\Xi$ from \eqref{E:formula_curly_B} is a non-decreasing function.

The load-free intervals are an integral part of the proof that $(\pmb{X}(q), q \ge 0)$ is a multiplicative coalescent.%
The advantage of working here with the modified reflected process $(\mathfrak{B}^{\pmb{x},q}, q \ge 0)$, is that we can couple it in a consistent way with a homogeneous Poisson point process 
$\Lambda$, as described below.

Let us first observe the following monotonicity property.
\begin{lemma}
\label{lem:curly_B_monotone}
    For any finite $\pmb{x}$ and for all $0\leq q_1<q_2$ we have
    $$
    \mathfrak{B}^{\pmb{x},q_1}(s) \leq \mathfrak{B}^{\pmb{x},q_2}(s), \ s \geq 0.
    $$
\end{lemma}
\begin{proof}
Recall the alternative expression \eqref{E:formula_curly_B} for $\mathfrak{B}^{\pmb{x},q}$.
In addition, since $
\Xi$ is non-decreasing, notice that $q \mapsto \Xi(\xi_i/q)$ is non-increasing in $q$, for each $i\ge 1$. 
This together with \eqref{E:formula_curly_B} clearly implies that
\(
    \mathfrak{B}^{\pmb{x},q_1}(s) \leq \mathfrak{B}^{\pmb{x},q_2}(s),
\)
for every $s \geq 0$.
\end{proof}

Let now $\Lambda$ be a homogeneous Poisson point process on $[0,\infty)^2$.
Consider the sets    
\begin{equation*}
\label{eq:def_mathfrak_A}
    \mathfrak{A}_i(q) = \{(x, y) \in \big[ \Xi(\xi_i/q) , \Xi(\xi_i/q) + X_i(q) \big) \times [0, \infty): y \le q \cdot \mathfrak{B}^{\pmb{x}, q}(x) \}.
\end{equation*}
Let $T_0 = 0$ and denote $T_1 < T_2 < \dots < T_{\mathrm{len}(\pmb{x}) - 1}$ the successive times of merging of the excursions of $(B^{\pmb{x},q}, \,q>0)$,  viewed as a process in $q$.
Note that, due to Lemma \ref{lem:curly_B_monotone}, if $q<q'$ then we have $\mathfrak{A}_i(q) \subset \mathfrak{A}_i(q')$ for all $i\geq 1$ in between two successive coalescence times (more precisely, almost surely on $\cup_{j\geq 0} \{T_j < q < q' < T_{j + 1}\}$).
Finally, let $(\pmb{X}(q), \pmb{\mathfrak{N}}(q))$ be the natural embedding of $$\big( (X_{1}(q), X_{2}(q), \dots), (\Lambda(\mathfrak{A}_{1}(q) ), \Lambda(\mathfrak{A}_{2}(q) ), \dots) \big)$$ into $\mathbb{U}_\searrow$, for every $q \ge 0$.
Recall that possible ties in excursion lengths are resolved using the corresponding surplus values.
\begin{lemma}
\label{lem:rate_increase_Lambda_along_B}
For any given $q>0$ and each $i$ (such that ${X}_i(q)>0$)
$$
\frac{\mathrm{d}}{ \mathrm{d} q} \operatorname{volume}(\mathfrak{A}_i(q)) = \frac{{X}_i(q)^2}{2}, \ \text{almost surely}.
$$
\end{lemma}

Let us prove Theorem \ref{thm:dynamicAMC2} before addressing Lemma \ref{lem:rate_increase_Lambda_along_B}.

\begin{proof}[Proof of Theorem \ref{thm:dynamicAMC2}]
Recall that $(\pmb{z}^{i, j}, \pmb{n}^{i,j})$  is obtained by merging components $i$ and $j$ into a component with volume $z_i + z_j$ and surplus $n_i + n_j$, and reordering the coordinates to obtain an element in $\mathbb{U}_\searrow$.
The process $(\pmb{X}(q), q \ge 0)$ evolves according to the multiplicative coalescent dynamics, see Proposition 5 in \cite{multcoalnew}.
According to the definition of $(\pmb{X}(q), q \ge 0)$, knowing that the process 
$(\pmb{X}, \pmb{\mathfrak{N}})$ is at state $(\pmb{z}, \pmb{n})$, it is clear that with rate $z_i \cdot z_j$, with $i \neq j$, the excursions of lengths $z_i$ and $z_j$ merge and the process jumps from $(\pmb{z}, \pmb{n})$ to $(\pmb{z}^{i, j}, \pmb{n}^{i,j})$, where reordering is included if necessary to obtain an element in $\mathbb{U}_\searrow$.

In addition, due to the definition $\Lambda(\mathfrak{A}_i)$ and Lemma \ref{lem:rate_increase_Lambda_along_B} we conclude that for each $i$
the transition from $(\pmb{z}, \pmb{n})$ to 
 $(\pmb{z}, \pmb{n}^{i})$ arrives at rate $z_i^2/2$, where $(\pmb{z}, \pmb{n}^{i})$ is the state obtained by changing only the $i^{\mathrm{th}}$ element in $\pmb{n}$ into $n_i+1$,
and reordering the coordinates if needed, to obtain an element in $\mathbb{U}_\searrow$. 
The vector $(\pmb{z}, \pmb{n}^{i})$ can be also defined as follows:
if $k(i) := \min\{j: (z_j, n_j)= (z_i,n_i)\}$ (note that $k(i) \le i$), then
$(\pmb{z}, \pmb{n}^{i})=\big( (z_j,n_j+\delta_{k(i)}(j))_{j \ge 1} \big)$, where $\delta_i(j)$ is the Kronecker delta indicator.
\end{proof}

\begin{proof}[Proof of Lemma \ref{lem:rate_increase_Lambda_along_B}]    
At this point, it would not be difficult to develop a proof of Lemma \ref{lem:rate_increase_Lambda_along_B} based on Lemma \ref{coro:area_below_B}.
However, we propose an alternative and self-contained method, independent from our study in the previous sections.

Let us compute $\omega$-by-$\omega$ the area below $q \cdot B^{\pmb{x}, q}$, restricted on a single excursion interval. 
Without loss of generality let us denote by $\{\pi_l, \pi_{l + 1}, \dots, \pi_r\}$ the set of vertices that form an excursion of length $X = x_{\pi_l} +  x_{\pi_{l + 1}} + \dots + x_{\pi_r}$.
Let us denote by $[\xi_{(l)}/q, \xi_{(l)}/q + X]$ this excursion interval.
Our goal is to compute the area below $B^{\pmb{x}, q}$, i.e.\ the integral 
\begin{equation}\label{eq:integralofB_during_excursion}
    \int_{\xi_{(l)}/q}^{\xi_{(l)}/q + X} B^{ \pmb{x}, q}(s) \mathrm{d}s = \int_{\xi_{(l)}/q}^{\xi_{(l)}/q + X} Z^{ \pmb{x}, q}(s) \mathrm{d}s - X \cdot Z^{ \pmb{x}, q}\left( \frac{\xi_{(l)}}{q}- \right),
\end{equation}   
where the second term on the right hand side is a consequence of the excursion property
\[
    \inf_{u \le s} Z^{\pmb{x}, q}(u) = Z^{ \pmb{x}, q}\left( \frac{\xi_{(l)}}{q} - \right), \; \forall s \in \left[ \frac{\xi_{(l)}}{q}, \frac{\xi_{(l)}}{q} + X \right].
\]
Furthermore, the definition of $Z^{\pmb{x},q}$ yields
\begin{equation}\label{eq:Zstar_excursion}
    Z^{ \pmb{x}, q}\left( \frac{\xi_{(l)}}{q} - \right) = \sum_{i = 1}^{l-1} x_{\pi_l} - \frac{\xi_{(l)}}{q}.
\end{equation}
Similarly, we can compute
\begin{align}
\int_{\xi_{(l)}/q}^{\xi_{(l)}/q + X} Z^{ \pmb{x}, q}(s) \mathrm{d}s & = \int_{\xi_{(l)}/q}^{\xi_{(l)}/q + X} \left(  \sum_{i = 1}^{\operatorname{len}(\pmb{x})} x_{\pi_i} \mathbb{1}_{\{ \xi_{(i)}/q \le s \}} - s\right)   \mathrm{d}s \nonumber \\
&= \sum_{i = 1}^{r} x_{\pi_i} \int_{\xi_{(l)}/q}^{\xi_{(l)}/q + X} \mathbb{1}_{\{ \xi_{(i)}/q \le s \}} \mathrm{d}s - \int_{\xi_{(l)}/q}^{\xi_{(l)}/q + X}  s \,  \mathrm{d}s \nonumber \\
    &= X \sum_{i = 1}^{l-1} x_{\pi_i} + \sum_{i = l}^r x_{\pi_i} \left( \frac{\xi_{(l)}}{q} + X - \frac{\xi_{(i)}}{q} \right) - \frac{1}{2} \left( \frac{\xi_{(l)}}{q} + X \right)^2 + \frac{1}{2} \left( \frac{\xi_{(l)}}{q}\right)^2
    \nonumber \\
    &= X \sum_{i = 1}^{l-1} x_{\pi_i} + \frac{1}{2} X^2 - \sum_{i = l}^r x_{\pi_i}  \frac{\xi_{(i)}}{q}. \label{eq:integralZ_during_excursion}
\end{align}
Inserting the RHS of \eqref{eq:Zstar_excursion} and \eqref{eq:integralZ_during_excursion} into \eqref{eq:integralofB_during_excursion} we get that
\begin{equation}\label{eq:integralofB_during_excursion2}
\int_{\xi_{(l)}/q}^{\xi_{(l)}/q + X} B^{ \pmb{x}, q}(s) \mathrm{d}s = \frac{1}{2} X^2 + X \frac{\xi_{(l)}}{q} - \sum_{i = l}^r x_{\pi_i} \frac{\xi_{(i)}}{q}.
\end{equation}
Recall the succesive times $(T_j)_j$ of merging of excursions of $B^{\pmb{x},\cdot}$.
For fixed $q > 0$ and $\epsilon > 0$, let us consider the event $C_\epsilon^q := \cup_{j = 1}^{\mathrm{len}(\pmb{x}) - 1} \{T_j < q - \epsilon < q + \epsilon < T_{j + 1}\}$.
Let $I_i$ be the random index such that the $i^{\text{th}}$ longest excursion of $B^{\pmb{x}, q}$ starts at $\xi_{(I_i)}/q$.
We can apply 
\eqref{eq:integralofB_during_excursion2} with $l=I_i$, 
 to conclude that on $C_\epsilon^q$
\begin{align*}
    \operatorname{volume}(\mathfrak{A}_i(q')) - \operatorname{volume}(\mathfrak{A}_i(q)) &= q' \int_{\xi_{(I_i)}/q'}^{\xi_{(I_i)}/q' + X_i} B^{\pmb{x}, q'} (s) \mathrm{d}s - q \int_{\xi_{(I_i)}/q}^{\xi_{(I_i)}/q + X_i} B^{\pmb{x}, q} (s) \mathrm{d}s \\ 
    &= \frac{1}{2} X_i^2 (q'-q), 
\end{align*}
for each $q'$ such that $|q' - q| < \epsilon$.
Since the event $\cup_{\epsilon > 0} C_\epsilon^q$ has probability $1$,
this concludes the proof.
\end{proof}

\section{Related results and open problems}
\label{S:Scaling}

Define the parameter space
\[
	\mathcal{I} := \big( (0, \infty) \times (- \infty, \infty) \times l^3_\searrow \big) \cup \big( \{0\} \times (- \infty, \infty) \times l^3_\searrow \setminus l^2_\searrow \big).
\]
Let $(\kappa, \tau, \pmb{c}) \in \mathcal{I}$. 
Given a family $(\xi_j')_j$ of independent exponential random variables, where $\xi_j'$ has rate $c_j$,
define 
\begin{equation*}
 V^{\pmb{c}} (s) = \sum_j \left(c_j \cdot \mathbb{1}_{(\xi_j' \leq s)} - c_j^2 \cdot s \right)
, \ s \geq  0.
\end{equation*}
For each $t\in \mathbb{R}$, let
$$
W^{\kappa,t- \tau, \pmb{c}}(s) = \kappa^{1/2}W(s) -\tau s - \frac{1}{2}\kappa s^2  + V^{\pmb{c}} (s) + t s, \ s \geq 0 \label{defWtc},
$$
where
$W$ is standard Brownian motion, and $W$ and $V^{\pmb{c}}$ are independent, and let
\begin{equation}\label{eq:def_B^kappa^tau^c}
 B^{\kappa, t-\tau,\pmb{c}}(s) := W^{\kappa, t-\tau,\pmb{c}}(s) - \min_{0 \leq s^\prime \leq s} W^{\kappa, t-\tau,\pmb{c}}(s^\prime), \ s \geq 0.
\end{equation}

For a given $\pmb{x} \in l^2_\searrow$ let
\(
	\sigma_r(\pmb{x}) := \sum_i x_i^r, \text{ for } r = 1,2,3.
\)
Assume that the sequence of initial masses $(\pmb{x}^{(n)})_n$ satisfies the following conditions:
\begin{align}
\frac{\sigma_{3}\left(\pmb{x}^{(n)}\right)}{\left(\sigma_{2}\left(\pmb{x}^{(n)}\right)\right)^{3}} &\rightarrow \kappa+\sum_{j} c_{j}^{3}, \label{eq:hypo1} \\
\frac{x_{j}^{(n)}}{\sigma_{2}\left(\pmb{x}^{(n)}\right)} & \rightarrow c_{j}, \quad j \geq 1, \label{eq:hypo2} \\
\sigma_{2}\left(\pmb{x}^{(n)}\right) & \rightarrow 0, \label{eq:hypo3}
\end{align}
as $n \to \infty$.
It is always possible to find sequences $\big( \pmb{x}^{(n)} \big)_{n \ge 1}$ satisfying \eqref{eq:hypo1} \eqref{eq:hypo2} and \eqref{eq:hypo3}, for every $(\kappa, \tau, \pmb{c}) \in \mathcal{I}$, see \cite[Lemma 8]{EBMC}.

For each $t$, let $\mathcal{X}(t)=\mathcal{X}^{\kappa,\tau,\pmb{c}}(t)$ be the infinite vector of ordered excursion lengths of $B^{\kappa, t-\tau, \pmb{c}}$ away from $0$. 
Theorem 1.2 in \cite{multcoalnew} says that the process $(\mathcal{X}(t), t\in (-\infty,\infty))$ is a realization of the extreme eternal
multiplicative coalescent\ corresponding to $(\kappa, \tau, \pmb{c}) \in \mathcal{I}$.
Furthermore, Corollary 11 in \cite{multcoalnew} says that $\mathcal{X}(t)$ is the scaling limit of the multiplicative coalescent started from $\pmb{x}^{(n)}$ (satisfying \eqref{eq:hypo1}, \eqref{eq:hypo2} and \eqref{eq:hypo3}), observed at time $q_n(t) = t + 1/\sigma_2(\pmb{x}^{(n)})$, when $n \to \infty$.

Our Theorems \ref{thm:encodingAMC} and \ref{thm:dynamicAMC2} and
Remarks \ref{rmk:joint_intensity}--\ref{rmk:joint_intensity2} are promising in view of novel scaling limits for near-critical random graphs, outside the domain of attraction of the Aldous standard multiplicative coalescent.
The eternal standard augmented multiplicative coalescent\ is the original one of Bhamidi et al.~\cite{bhamidietal2}, a version of which was constructed in \cite{bromar15} as the scaling limit of the random graph with surplus counts for special initial configurations of the form
\begin{equation*}
    x_1 = x_2 = \ldots = x_n = 1/n^{2/3}, \text{ and } 0 = x_{n+1} = \ldots, \text{ as } n \text{ diverges.} 
\end{equation*}
Also, Dhara et al.\ \cite[Thm.\ 3.6]{Dhara2017} obtained a version of the eternal standard augmented multiplicative coalescent for configuration models with finite third moment degrees.
Furthermore, \cite[Thm.\ 2 and 5]{Dhara2019} found the scaling limit for the sizes of the connected components and the number of surplus edges for configuration models with heavy-tailed distribution in the degree.
In this latter setting, the version of the augmented multiplicative coalescent that appears is not the standard one, instead it is the one where the Brownian part vanishes, (i.e.\  $\kappa = 0$ and $\pmb{c} \in l^3_\searrow \setminus l^2_\searrow$).
At present, as far as we know, these are the only versions of the augmented multiplicative coalescent that have been studied. 

Let $\Lambda$ be a homogeneous Poisson point process on $[0,\infty)\times [0,\infty)$, independent of $\sigma\{W,V^{\pmb{c}}\}$.
In analogy to \cite{aldRGMC,bhamidietal2}, let $\Lambda^{\kappa, t-\tau,\pmb{c}}(s)$ be the number of points in $\mathcal{N}$
 below the curve $u\mapsto B^{\kappa,t-\tau,\pmb{c}}(u)$, $u\in [0,s]$. 
To each excursion of $B^{\kappa, t-\tau, \pmb{c}}$ above $0$, one can assign a random ``mark count'' as $\Lambda^{\kappa, t-\tau,\pmb{c}}$  of the region below the curve $B^{\kappa, t-\tau, \pmb{c}}$ restricted to this excursion (see \cite[Section 2.3.2]{bhamidietal2} for details in the standard setting). 
Let $\mathcal{N}_i(t)$ be this count corresponding to the $i^{\text{th}}$ longest excursion of $B^{\kappa,t-\tau,\pmb{c}}$, and define $\mathcal{N}(t) = \mathcal{N}^{\kappa,\tau,\pmb{c}}(t) =(\mathcal{N}_1(t),\mathcal{N}_2(t),\ldots)$. 
Given the observations made in previous sections, the following can be anticipated 
(a work in progress by the authors is devoted the proof of this claim):
\begin{conjecture}
\label{Thm:AMC_scaling_limit}
Fix a $(\kappa,\tau,\pmb{c}) \in \mathcal{I}$. Then $\big( (\mathcal{X}(t),\mathcal{Y}(t))$, $ \,-\infty < t< \infty) \big)$ is a càdlàg realization of the {\em eternal  augmented multiplicative coalescent} corresponding to $(\kappa,\tau,\pmb{c})$. 
Furthermore, $\big( (\mathcal{X}(t),\mathcal{Y}(t))$, $ \,-\infty < t< \infty \big)$ is the simultaneous scaling limit of near-critical random graph component sizes and surplus counts,  under the hypotheses of the initial configurations \eqref{eq:hypo1}, \eqref{eq:hypo2} and \eqref{eq:hypo3}. 

In addition, the extreme eternal augmented multiplicative coalescents are only the constant ones, and the non-trivial ones given here (corresponding to valid parameters $(\kappa,\tau,\pmb{c})$). Any eternal augmented multiplicative coalescent\ is a mixture of extreme ones.
\end{conjecture}

The excursion mosaic and the accompanying PPP family $\zeta^{l;j-k}$, $j\leq k\leq l$ (see Section \ref{S:surtopforo}) has a much richer structure than the mere component sizes superimposed by surplus edge counts. 
Is there a natural framework and candidate for its scaling limit in the near-critical regime(s)? 
This insight would surely encompass a clearer understanding of mark counts $\mathcal{N}$  in the eternal augmented coalescent.


\providecommand{\bysame}{\leavevmode\hbox to3em{\hrulefill}\thinspace}
\providecommand{\MR}{\relax\ifhmode\unskip\space\fi MR }
\providecommand{\MRhref}[2]{%
	\href{http://www.ams.org/mathscinet-getitem?mr=#1}{#2}
}
\providecommand{\href}[2]{#2}

\end{document}